\documentclass{article}
\usepackage{amssymb}
\usepackage{amsmath,color}
\usepackage[hidelinks]{hyperref}

\usepackage[square,numbers]{natbib}         
\usepackage{ntheorem}
\newtheorem{lem}{Lemma}[section]
\newtheorem{cor}[lem]{Corollary}
\newtheorem{teo}[lem]{Theorem}

\newtheorem{os}[lem]{Remark}
\newtheorem{defi}[lem]{Definition}
\newtheorem{prop}[lem]{Proposition}

\newcommand{\qed}{\thinspace\null\nobreak\hfill\hbox{\vbox{\kern-.2pt\hrule
			height.2pt depth.2pt\kern-.2pt\kern-.2pt \hbox to2.5mm{\kern-.2pt\vrule
				width.4pt \kern-.2pt\raise2.5mm\vbox to.2pt{}\lower0pt\vtop
				to.2pt{}\hfil\kern-.2pt \vrule
				width.4pt \kern-.2pt}\kern-.2pt\kern-.2pt\hrule height.2pt depth.2pt
			\kern-.2pt}}\par\medbreak}

\newcommand{\R}{\mathbb{R}}

\newcommand{\C}{\mathbb{C}}

\newcommand{\N}{\mathbb{N}}

\newcommand{\Rp}{\textrm{\emph{Re}\,}}
\newcommand{\Ip}{\textrm{\emph{Im}\,}}
\newcommand{\eps}{\varepsilon}

\newcommand{\ds}{\displaystyle}

\date{}
\begin{document}

\title{
Kernel bounds for parabolic operators having first-order degeneracy at the boundary}
\author{L. Negro \thanks{Dipartimento di Matematica e Fisica ``Ennio De Giorgi'', Universit\`a del Salento, C.P.193, 73100, Lecce, Italy.
e-mail:  luigi.negro@unisalento.it} \qquad C. Spina \thanks{Dipartimento di Matematica e Fisica ``Ennio De Giorgi'', Universit\`a del Salento, C.P.193, 73100, Lecce, Italy.
e-mail:  chiara.spina@unisalento.it}}

\maketitle
\begin{abstract}
\noindent 
We study kernel estimates for parabolic problems governed by singular elliptic   operators 
 \begin{equation*}
\sum_{i,j=1}^{N+1}q_{ij}D_{ij}+c\frac{D_y}{y},\qquad c+1>0, 
\end{equation*} 
in the half-space $\R^{N+1}_+=\{(x,y): x \in \R^N, y>0\}$ under Neumann boundary conditions at $y=0$.  

\par

\end{abstract}

\section{Introduction}
In this paper   we  study upper  estimates for the parabolic heat kernel of the   degenerate   operator 
 \begin{equation}\label{def L}
	\mathcal L =\sum_{i,j=1}^{N+1}q_{ij}D_{ij}+c\frac{D_y}{y}
\end{equation}
in the half-space $\R^{N+1}_+=\{(x,y): x \in \R^N, y>0\}$, under Neumann boundary conditions at $y=0$.  Here  $Q=\left(q_{ij}\right)_{i,j=1,\dots,N+1}\in\R^{N+1,N+1}$ is a constant and positive definite  matrix and $c\in\R$ satisfies $\frac c\gamma +1>0$ where $\gamma=q_{N+1,N+1}$.

\medskip
In the special case $$\mathcal L=\Delta_x+B_y, \qquad  B_y=D_{yy}+ \frac cy D_y,$$  that is when the mixed second order derivatives do not appear in the operator, the  one-dimensional Bessel operator $B_y$ and the Laplacian $\Delta_x$ commute and sharp kernel estimates  can be  proved  after an explicit description of the kernel of the Bessel operator (see \cite{MNS-Caffarelli}). These operators  play a major role in the investigation of the fractional powers $(-\Delta_x)^s$ and  $(D_t-\Delta_x)^s$, $s=(1-c)/2$, through the  ``extension procedure" of Caffarelli and Silvestre, see \cite{Caffarelli-Silvestre}.  The strategy above based on the commutativity of $\Delta_x$ and $B_y$ clearly fails for the general operator \eqref{def L} due to the presence of the mixed derivatives $D_{xy}$. On the other hand,  adding this term  is   crucial  for treating degenerate operators in domains, through a localization procedure.
\medskip 
	
	For this reason  we investigate, in this work, the validity of upper kernel estimate for the  general operator $\mathcal L$ defined in \eqref{def L} by  revisiting the classical method based upon the equivalence between Gagliardo-Nirenberg type inequalities and ultracontractivity estimates.

\medskip
 Before describing the  strategy here employed, we remark that elliptic and parabolic solvability of the associated problems in weighted $L^p$ spaces have been already investigated in 
\cite{MNS-Singular-Half-Space, MNS-Degenerate-Half-Space, MNS-Caffarelli, MNS-CompleteDegenerate} where we proved that, under suitable assumptions on $m$ and $p$ and for $c>-1$,  $\mathcal L$, endowed with  Neumann boundary conditions, generates an analytic semigroup in $L^p_m$ $=L^p(\R^{N+1}_+; y^m dxdy)$. In addition  we characterized its domain as a weighted Sobolev space and showed that it has the  maximal regularity.  In addition we also refer to   \cite{dong2020parabolic, dong2021weighted, dong2020RMI, dong2020neumann} where the authors studied  operators of this form, even for variable coefficients, using tools form linear PDE and Muckhenoupt weights. In \cite{Robinson-Sikora2008}, among other properties, kernel estimates for some second order degenerate operators in divergence form are considered.
In our knowledge, kernel estimates for the operators in the general form \eqref{def L} are new. 


\medskip
Our main results are Theorems \ref{complexGauss} and \ref{complete-oblique} where, if  $\frac c \gamma+1>0$, we prove that $\mathcal L$ generates a contractive analytic semigroup  in $L^2_{\frac c\gamma}$	and its heat kernel $p_{{\mathcal L}}$, written  with respect the measure $y^\frac{c}{\gamma}dz$, satisfies the upper  Gaussian estimate
\begin{align*}
	0\leq p_{{\mathcal  L}}(t,z_1,z_2)
	\leq C t^{-\frac{N+1}{2}} y_1^{-\frac{c}{2\gamma}} \left(1\wedge \frac {y_1}{\sqrt t}\right)^{\frac{c}{2\gamma}} y_2^{-\frac{c}{2\gamma}} \left(1\wedge \frac{y_2}{\sqrt t}\right)^{\frac{c}{2\gamma}}\,\exp\left(-\dfrac{|z_1-z_2|^2}{kt}\right),
\end{align*}
where $t>0$, $z_1=(x_1,y_1),\ z_2=(x_2,y_2)\in\R^{N+1}_+$;  we also extend the estimate to complex times.
As consequence of our kernel  estimates, we deduce, after a linear transformation, upper kernel estimates also for the operator 
\begin{equation}\label{def L oblique}
	\mathcal L =\sum_{i,j=1}^{N+1}q_{ij}D_{ij}+\frac{b\cdot\nabla_x}{y}+\frac{cD_y}{y},\qquad b\in\R^N,\quad c\neq 0,
\end{equation}
 under an oblique  boundary conditions at $y=0$.

\medskip
 We now describe the  strategy of proof employed. A suitable linear change of variables, allows to deal with the  equivalent  degenerate   operator 
\begin{equation*}  
	\mathcal L :=
	\Delta_x u+2a\cdot \nabla_xD_yu+ B_yu ,\qquad a=(a_1, \dots, a_N) \in\R^N,\ |a|<1.
\end{equation*}
Note that the condition $|a| <1$ is equivalent to the ellipticity of the top order coefficients. 

\smallskip
We start by constructing, in Section \ref{L2},  using form methods, the semigroup generated by $\mathcal L$ in a weighted $L^2$-space, where the weighted measure  $y^c\,dx\,dy$ takes into account the degeneracy of the operator. We prove continuity, accretivity and sectoriality for the form associated with $\mathcal L$ and its adjoint. 

Then, we try to use the standard machinery, relying on Gagliardo-Nirenberg
type inequalities, to deduce the ultracontractivity of the semigroup. To this aim, we investigate in Section \ref{Section weighted sobolev} the validity of some anisotropic  weighted Sobolev inequalities. Some difficulties appear  since the semigroup $e^{t\mathcal L}$ is not ultracontractive when $c<0$: this is due to the fact that the   Gagliardo-Nirenberg inequalities  are valid only locally due to the singularity of the operator at the boundary  $y=0$.

To overcome this obstruction, we introduce, in Section \ref{Section Kernel}, a further  weight $\phi$ 
which allow us to take into account the degeneracy of the operator at the boundary. This introduction is essential for $c<0$ since allows to prove  Gagliardo-Nirenberg inequalities in the whole half-space for any $c>-1$. When $c\geq 0$, this approach, although not mandatory,  is needed in order to deduce precise heat kernel estimates (see Remark \ref{Oss c>0}). Another obstruction to the application of the standard machinery  is  the non-symmetry of the operator $\mathcal L$: this force to study at the same time also the adjoint operator $\mathcal L^\ast$, proving $L^2-L^\infty$ estimate for $e^{t\mathcal L}$ and $e^{t\mathcal L^\ast}$. Then we implement Davies' trick to add the  Gaussian factor. Moreover we extend the estimates to complex times. We do not prove lower bounds; however, in the case $\Delta_x+B_y$,  our upper bounds coincide with those of \cite{MNS-Caffarelli} which are  optimal.
 
 Finally in Section \ref{Section oblique} we deduce upper kernel estimates  for the operator in  \eqref{def L oblique} having oblique boundary condition.

\bigskip
\noindent\textbf{Notation.} For $N \ge 0$, $\R^{N+1}_+=\{(x,y): x \in \R^N, y>0\}$. We write $\nabla u, D^2 u$ for the gradient and the Hessian matrix of a function $u$ with respect to all $x,y$ variables and $\nabla_x u, D_y u, D_{x_ix_j }u, D_{x_i y} u$ and so on, to distinguish the role of $x$ and $y$.

For $m \in \R$ we consider the measure $y^m dx dy $ in $\R^{N+1}_+$ and  we write $L^p_m(\R_+^{N+1})$, and often only $L^p_m$ when $\R^{N+1}_+$ is understood, for  $L^p(\R_+^{N+1}; y^m dx dy)$. 


$\C^+$ stands for $\{ \lambda \in \C: \Rp \lambda >0 \}$ and, for $|\theta| \leq \pi$, we denote by  $\Sigma_{\theta}$  the open sector $\{\lambda \in \C: \lambda \neq 0, \ |Arg (\lambda)| <\theta\}$.  

Given $a$ and $b$ $\in\R$, $a\wedge b$, $a \vee b$  denote  their minimum and  maximum. We  write $f(x)\simeq g(x)$ for $x$ in a set $I$ and positive $f,g$, if for some $C_1,C_2>0$ 
\begin{equation*}
	C_1\,g(x)\leq f(x)\leq C_2\, g(x),\quad x\in I.
\end{equation*}

\bigskip
\noindent\textbf{Acknowledgment.}
The authors are members of the INDAM (``Istituto Nazionale di Alta Matematica'') research group GNAMPA (``Gruppo Nazionale per l’Analisi Matematica, la Probabilità
e le loro Applicazioni'').

\section{The operator $\mathcal L=\Delta_{x} +2a\cdot\nabla_xD_y+ B_yu$ in $L^2_{c}$} \label{L2}

Let $c\in\R$, $a=(a_1, \dots, a_N) \in\R^N$ such that $c+1>0$, $|a|<1$. In this section we study the $L^2$ theory related to  the  degenerate   operator 
\begin{equation}  \label{La}
	\mathcal L :=\Delta_{x} +2\sum_{i=1}^Na_{i}D_{x_iy}+D_{yy}+\frac{c}{y}D_y=\Delta_x u+2a\cdot \nabla_xD_yu+ B_yu 
\end{equation}
in $L^2_c$ equipped with Neumann boundary condition.  The requirement $|a| <1$ is equivalent to the ellipticity of the top order coefficients.

We use the Sobolev space $H^{1}_{c}:=\{u \in L^2_{c} : \nabla u \in L^2_{c}\}$ equipped with the inner product
\begin{align*}
	\left\langle u, v\right\rangle_{H^1_{c}}:= \left\langle u, v\right\rangle_{L^2_{c}}+\left\langle \nabla u, \nabla v\right\rangle_{L^2_c}.
\end{align*}
%
The condition $c+1>0$ assures, by \cite[Theorem 4.9]{MNS-Sobolev}, that the set 
\begin{equation} \label{defC}
	\mathcal{C}:=\left \{u \in C_c^\infty \left(\R^N\times[0, \infty)\right), \ D_y u(x,y)=0\  {\rm for} \ y \leq \delta\ {\rm  and \ some\ } \delta>0\right \},
\end{equation}
is dense in $H^1_c$. Moreover,  by \cite[Remark 4.14]{MNS-Sobolev}, if a function $u\in H^1_c$ has  support in $\R^N\times[0,b]$, then there exists a sequence $\left(u_n\right)_{n\in\N}\in\mathcal C$  such that $ \mbox{supp }u_n\subseteq \R^N\times[0,b]$  and  $u_n\to u$ in $H^1_c$.

We consider the  form  in $L^2_{c}$ 
\begin{align*}
	\mathfrak{a}(u,v)
	&:=
	\int_{\R^{N+1}_+} \langle \nabla u, \nabla \overline{v}\rangle\,y^{c} dx\,dy+2\int_{\R^{N+1}_+} D_yu\, a\cdot\nabla_x\overline{v}\,y^{c} dx\,dy , \quad
	D(\mathfrak{a})=H^1_{c}
\end{align*}

and its adjoint $\mathfrak{a}^*(u,v)=\overline{\mathfrak{a}(v,u)}$  
\begin{align*}
	\mathfrak{a^*}(u,v)=\overline{\mathfrak{a}(v,u)}
	&:=
	\int_{\R^{N+1}_+} \langle \nabla u, \nabla \overline{v}\rangle\,y^{c} dx\,dy+2\int_{\R^{N+1}_+} a\cdot \nabla_x u\,D_y\overline{v}\,y^{c} dx\,dy.
\end{align*}

\begin{prop} \label{prop-form}
	The forms $\mathfrak{a}$, $\mathfrak{a^*}$ are continuous, accretive and   sectorial.
\end{prop}
{\sc Proof.} We consider only the form $\mathfrak{a}$ since the adjoint form can be handled similarly. If $u\in H^1_{c}$
\begin{align*}
	\Rp \mathfrak{a}(u,u)\geq \|\nabla_x u\|^2_{L^2_c}+\|D_y u\|^2_{L^2_c}-2|a| \|\nabla_x u\|_{L^2_c}\|D_y u\|_{L^2_c}\geq (1-|a|)(\|\nabla_x u\|^2_{L^2_c}+\|D_y u\|^2_{L^2_c}).
\end{align*}
By the ellipticity assumption $|a|<1$, the accretivity follows.
Moreover 
\begin{align*}
	|\Ip \mathfrak{a}(u,u)|\leq  2|a|\|\nabla_x u\|_{L^2_c}\|D_y u\|_{L^2_c}\leq |a|(\|\nabla_x u\|^2_{L^2_c}+\|D_y u\|^2_{L^2_c})\leq \frac{|a|}{(1-|a|)}\Rp \mathfrak{a}(u,u).
\end{align*}
This proves the sectoriality and then the continuity of the form. \qed

We define the operators $\mathcal L$ and $\mathcal L^*$ associated respectively to the forms $\mathfrak{a}$ and $\mathfrak{a}^*$  by
\begin{align} \label{BesselN}
	\nonumber D( \mathcal L)&=\{u \in H^1_{c}: \exists  f \in L^2_{c} \ {\rm such\ that}\  \mathfrak{a}(u,v)=\int_{\R^{N+1}_+} f \overline{v}y^{c}\, dz\ {\rm for\ every}\ v\in H^1_{c}\},\\  \mathcal Lu&=-f;
\end{align}
\begin{align} \label{adjoint}
	\nonumber D( \mathcal L^*)&=\{u \in H^1_{c}: \exists  f \in L^2_{c} \ {\rm such\ that}\  \mathfrak{a}^*(u,v)=\int_{\R^{N+1}_+} f \overline{v}y^{c}\, dz\ {\rm for\ every}\ v\in H^1_{c}\},\\  \mathcal L^*u&=-f.
\end{align}
If $u,v$ are smooth function with compact support in the closure of $\R_+^{N+1}$ (so that they do not need to vanish on the boundary), it is easy to see integrating by parts that $$-\mathfrak a (u,v)= \langle \Delta_x u+2a\cdot \nabla_xD_yu+ B_yu, \overline v\rangle_{L^2_c}  $$
if $\lim_{y \to 0} y^c D_yu(x,y)=0$. This means that 
$\mathcal L$ is the operator $\Delta_x +2a\cdot \nabla_xD_y+ B_y$ with Neumann boundary conditions at $y=0$. On the other hand 
$$-\mathfrak a^* (u,v)= \left\langle \Delta_x u+2a\cdot \nabla_xD_yu+2c\frac {a\cdot \nabla_x u}{y}+ B_yu, \overline v\right\rangle_{L^2_c}  $$
if $\lim_{y \to 0} y^c \left ( D_yu(x,y)+2a\cdot \nabla_x u(x,y)\right )=0$ and therefore $\mathcal L^*$ is the operator $\Delta_x +2a\cdot \nabla_xD_y+2c\frac {a\cdot \nabla_x u}{y}+ B_y$ with the above oblique condition at $y=0$.  We refer the reader to \cite{MNS-Sobolev} for further details about the Sobolev spaces $H^1_c$ and their boundary conditions at $y=0$.
\begin{prop}\label{generation L2}
	$\mathcal L$ and $\mathcal L^*$  generate contractive analytic semigroups $e^{z \mathcal L}$, $e^{z \mathcal L^*}$, $z\in\Sigma_{\frac{\pi}{2}-\arctan \frac{|a|}{1-|a|}}$,  in $L^2_{c}$. Moreover the semigroups $(e^{t\mathcal L})_{t \geq 0}, (e^{t\mathcal L^*})_{t \geq 0}$ are positive and $L^p_c$-contractive for $1 \leq p \leq \infty$.
\end{prop}
{\sc Proof.} We argue only for $\mathcal L$.
The  generation  result immediately follows from   Proposition \ref{prop-form} and \cite[Theorem 1.52]{Ouhabaz}.
The positivity follows by \cite[Theorem 2.6]{Ouhabaz} after observing that, if $u\in H^1_{c}$, $u$ real, then $u^+\in H^1_{c}$ and
\begin{align*}
	\mathfrak{a}(u^+,u^-)
	&:=
	\int_{\R^{N+1}_+} \langle \nabla u^+, \nabla u^-\rangle\,y^{c} dx\,dy+2\int_{\R^{N+1}_+}D_yu^+ a\cdot \nabla_xu^-\,y^{c} dx\,dy=0.
\end{align*}
Finally, the $L^\infty$-contractivity follows by \cite[Corolary 2.17]{Ouhabaz}  after observing that if $0\leq u\in H^1_{c}$,  then $1\wedge u, (u-1)^+\in H^1_{c}$ and, since $\nabla (1\wedge u)=\chi_{\{u<1\}}\nabla u$ 
	and 
	$\nabla (u-1)^+=\chi_{\{u>1\}}\nabla u$, one has 
	\begin{align*}
		\mathfrak{a}(1\wedge u,(u-1)^+)&=0.
	\end{align*}
	\qed

The previous proposition can be generalized by showing that the above semigroup extrapolates to the weighted spaces $L^p\left(\R^{N+1}_+, y^mdxdy\right)$, $m\in\R$, when   $0<\frac{m+1}p<c+1$. We refer the reader to \cite{MNS-Singular-Half-Space, MNS-Degenerate-Half-Space} where a characterization of the domain of the generators is given in terms of  suitable anisotropic weighted Sobolev spaces (see \cite{MNS-Sobolev}) and where it is also proved that the generated semigroups satisfy parabolic maximal regularity.

\medskip

Finally we collect in the following proposition the scaling properties of $\mathcal L$ and $e^{t\mathcal L}$ whose proofs follow after an easy straightforward computation.

\begin{prop}\label{Prop scaling} The scale homogeneity of $\mathcal L$ is  $2$:
		\begin{align}\label{scale.op}
			s^2\mathcal L=I_s^{-1}\mathcal L I_s, \qquad I_su(x,y)=u(sx,sy), \qquad s>0. 
		\end{align}
	Moreover the semigroup $e^{t\mathcal L}$ generated by $\mathcal L$  in $L^2_{c}$
		satisfies
		\begin{equation}\label{scale.semi}
			e^{s^2 z\mathcal L}=I_{s}^{-1}e^{-z\mathcal L}I_s, \qquad   z\in\Sigma_{\frac{\pi}{2}-\arctan \frac{|a|}{1-|a|}},\quad s>0. 
		\end{equation}
The same results  hold for the adjoint operator $\mathcal L^\ast$ and its generated semigroup. 
\end{prop}

\section{Weighted Sobolev  inequalities}\label{Section weighted sobolev}
In this section we investigate the validity of anisotropic weighted Sobolev inequalities for the space $H^1_{c}(\R^{N+1}_+)$. An essential role is played by the sign of $c$.  Indeed assuming that for some $p>2$ one has 
\begin{align*}
	H^{1}_{c}(\R^{N+1}_+)\hookrightarrow L^{p}_c\left(\R^{N+1}_+\right),\qquad 	\| u\|_{L^{p}_{c}\left(\R^{N+1}_+\right)}
	\leq C\|\nabla u\|_{L^2_{c}\left(\R^{N+1}_+\right)}.
\end{align*} 
then a simple scaling argument implies that the right expression of the Sobolev exponent $p$ is 
\begin{align*}
	2^{\ast}_{c}\stackrel{def}=2\frac{N+1+c}{N+c-1},\qquad\text{i.e.}\qquad \frac 1 {2^{*}_c}=\frac 1 2 -\frac{1}{N+1+c}.
\end{align*} 
Then a necessary condition for the validity of the latter inequalities is  $c\geq 0$: this assumption guarantees that  $2^\ast_c$ is not greater than the Sobolev exponent $2^{\ast}=2\frac{N+1}{N-1}$ for the Lebesgue measure since otherwise  this  would imply better Sobolev embedding on sets far from the singularity $y=0$, which is impossible.

For this reason we distinguish between the cases $c\geq 0$ and $c<0$. When $c\geq 0$ we prove that  the  embedding $H^1_{c}(\R^{N+1}_+)\hookrightarrow L_c^q(\R^{N+1}_+)$ holds for any $2\leq q<{2^{*}_c}$. On the other hand, when $c<0$, we prove that the Sobolev immersion  holds  only locally and with the same Sobolev exponent $2^\ast$ as for the Lebesgue measure.

\medskip 
We start by the  case $c\geq 0$. We need some preparation. We treat first the case $c=0$ where the latter embedding is just  the standard Sobolev inequality which we state below in the form we use.

\begin{lem}\label{GN}
	Let $N\in\N$. Then 
	for every $q\in (2,\infty)$ satisfying $\frac{1}{q}\geq \frac{1}{2}-\frac{1}{N+1}$,
	there exists $\widetilde{C}_q>0$ such that 
	for every $u\in H^1(\R^{N+1}_+)$,
	\begin{align*}
		\|u\|_{L^q(\R^{N+1}_+)}\leq 
		\widetilde{C}_{q}
		\|\nabla v\|_{L^2(\R^{N+1}_+)}^{(N+1)\left (\frac{1}{2}-\frac{1}{q} \right )}\|v\|_{L^2(\R^{N+1}_+)}^{1-(N+1)\left (\frac{1}{2}-\frac{1}{q} \right )}. 
	\end{align*}
\end{lem}

We state now  the following lemmas which can  be found in \cite{mazya}.

\begin{lem}  \label{Mazya 1}
	If $\alpha\in\R$ and  $q\geq 1$ satisfies
	\begin{align*}
		1\leq q\leq (N+1)'=\frac{N+1}{N},\qquad \beta=\alpha-1+\frac{q-1}{q}(N+1)>-\frac{1}{q}
	\end{align*}
	(in particular $\alpha+1>0$) then there exists a positive constant $C$ such that
	\begin{align}\label{Mazya 1 eq}
		\|y^{\beta} u\|_{L^q(\R^{N+1}_+)}\leq C\|y^{\alpha}\nabla u\|_{L^1(\R^{N+1}_+)},\quad \forall u\in C_c^\infty\left(\overline{\R^{N+1}_+}\right).
	\end{align}
	
\end{lem}
{\sc{Proof.}} See  \cite[Corollary 1, Section 2.1.7, page 138]{mazya}
\qed

\begin{lem}\label{Mazya 2}
	Let  $q, \ \beta,\  \alpha\in \R$, satisfying the following assumptions:
	\begin{itemize}
		\item[(i)] $2\leq q\leq 2^{\ast}= 2\frac{N+1}{N-1}$ \quad ($2\leq q$ when $N=1$);
		\item[(ii)] $\beta+\frac 1 q>0$;
		\item[(iii)] $\beta= \alpha-1+(N+1)\left(\frac 1 2 -\frac 1 q\right)$ \quad i.e. \quad $\frac 1 q=\frac 1 2- \frac{\beta-\alpha+1}{N+1}$.
	\end{itemize} 
	Then one has  for some positive constant $C$
	\begin{align*}
		\| y^\beta u\|_{L^q\left(\R^{N+1}_+\right)}
		\leq C\|y^{ \alpha}\nabla u\|_{L^2\left(\R^{N+1}_+\right)},\qquad \text{for all }\quad u\in C_c^\infty\left(\overline{\R^{N+1}_+}\right).
	\end{align*}
\end{lem}
{\sc{Proof.}} 
The proof follows by applying Lemma \ref{Mazya 1} with $q$, $\alpha$, $\beta$  replaced  respectively by $\frac{2q}{q+2}$, $\alpha+\frac{q\beta}2$, $\left(\frac q2+1\right)\beta$ to the function $|u|^{\frac{q}2+1}$ which yields
\begin{align*}
	\|y^{\beta} u\|^{\frac q 2+1}_{L^q(\R^{N+1}_+)}\leq C\frac{q+2}{2}\|y^{\alpha}|\nabla u|y^{\frac{q\beta}2}|u|^{\frac q 2}\|_{L^1(\R^{N+1}_+)}.
\end{align*}
Applying then the  Cauchy-Schwarz inequality to the right hand side yields
\begin{align*}
	\|y^{\beta} u\|^{\frac q 2+1}_{L^q(\R^{N+1}_+)}\leq C\frac{q+2}{2}\|y^\alpha\nabla u\|_{L^2(\R^{N+1}_+)}\|y^\beta u\|^{\frac q 2}_{L^q(\R^{N+1}_+)}
\end{align*}
which is the required claim.

\qed

We can now prove the    Sobolev immersion 
$H^{1}_{c}(\R^{N+1}_+)\hookrightarrow L^{2^\ast_c}_c\left(\R^{N+1}_+\right)$ assuming 
\begin{align*}
	2^{\ast}_{c}=2\frac{N+1+c}{N+c-1},\qquad c\geq 0.
\end{align*}

\begin{prop}\label{Sobolev1}
	Let  $c>0$ if $N=1$, $c\geq 0$ otherwise. Then   for some positive constant $C=C(c,N)$ one has
	\begin{align*}
		\| u\|_{L^{2^\ast_{c}}_{c}\left(\R^{N+1}_+\right)}
		\leq C\|\nabla u\|_{L^2_{c}\left(\R^{N+1}_+\right)},\qquad \text{for all }\quad u\in H^1_{c}\left(\R^{N+1}_+\right).
	\end{align*}  
\end{prop}
{\sc{Proof.}} Since the set $\mathcal C$ defined in \eqref{defC} is dense in $H^1_c$, we can assume without any loss of generality $u\in\mathcal C$.  The proof then  follows using Lemma   \ref{Mazya 2} with  $ \alpha=\frac{c}2$, $\beta=\frac {c} q$.\\
\qed

The following Gagliardo-Nirenberg type inequalities easily follow.

\begin{prop}\label{CKN2}
	Let $c\geq 0$  and let us fix  $q\in (2,\infty)$ satisfying  satisfying $\frac{1}{q}\geq \frac{1}{2}-\frac{1}{N+1+c}$  i.e.
	\begin{itemize}
		\item[(i)] $2< q\leq 2^\ast_c$ \quad if  $N>1$ and $c\geq 0$,   or $N=1$ and $c>0$;
		\item[(ii)] $q>2 $  \quad  if  $N=1$ and $c=0$.
	\end{itemize}
	Then  there exists a  positive constant $C$, which may depends on $c, N$ (and also on $q$ in case (ii)), such that 
	\begin{align*}
		\|u\|_{L^q_c(\R^{N+1}_+)}\leq C \| \nabla u\|_{L^2_{c}(\R^{N+1}_+)}^{(N+1+c)\left(\frac{1}{2}-\frac{1}{q}\right)}\| u\|_{L^2_c(\R^{N+1}_+)}^{1-(N+1+c)\left(\frac{1}{2}-\frac{1}{q}\right)},\qquad \text{for all }\quad u\in H^1_{c}\left(\R^{N+1}_+\right).
	\end{align*}
\end{prop}
{\sc Proof.} Case (ii) is just Lemma \ref{GN}. To prove (i) let us fix $u\in H^1_{c}\left(\R^{N+1}_+\right)$. Then, by Holder's inequality and by Proposition \ref{Sobolev1}, one has for some $C(c,N)$
\begin{align*}
	\|u\|_{L^q_c(\R^{N+1}_+)}&\leq \| u\|_{L^{2_c^\ast}_{c}(\R^{N+1}_+)}^{(N+1+c)\left(\frac{1}{2}-\frac{1}{q}\right)}\| u\|_{L^2_c(\R^{N+1}_+)}^{1-(N+1+c)\left(\frac{1}{2}-\frac{1}{q}\right)}\\[1ex]
	&\leq \Big[C\| \nabla u\|_{L^{2}_{c}(\R^{N+1}_+)}\Big]^{(N+1+c)\left(\frac{1}{2}-\frac{1}{q}\right)}\| u\|_{L^2_c(\R^{N+1}_+)}^{1-(N+1+c)\left(\frac{1}{2}-\frac{1}{q}\right)}
\end{align*}
which is the required estimates.\\\qed

\bigskip 
When  $c<0$, as explained before,  the immersion $H^1_{c}(\R^{N+1}_+)\hookrightarrow L^{2^\ast_c}_c\left(\R^{N+1}_+\right)$ fails since otherwise would imply better Sobolev embedding on sets far from the singularity $y=0$. 
In the next lemma we show that, in this case, the Sobolev immersion still holds locally (note that in this case  the Sobolev exponent for $|y|^c \, dz$ is the same as for the Lebesgue measure). 


\begin{lem}\label{Sobolev2}
	Let    $c+1>0$. The following properties hold.
	\begin{itemize}
		\item[(i)]If $N>1$ then  one has for some positive constant $C$
		\begin{align*}
			\|u\|_{L^{2^\ast}_{c}(\R^{N+1}_+)}\leq C\, \|y^{-\frac{c}{N+1}}\nabla u\|_{L^2_{c}(\R^{N+1}_+)},\qquad \text{for all }\quad u\in C_c^\infty\left(\overline{\R^{N+1}_+}\right).
		\end{align*}
		\item[(ii)]If $N=1$ then for every $q\geq 2$ there exists a  positive constant $C$ such that  \begin{align*}
			\|u\|_{L^{q}_{c}(\R^{N+1}_+)}\leq C\, \|y^{\frac{c+2}q-\frac c 2}\nabla u\|_{L^2_{c}(\R^{N+1}_+)},\qquad \text{for all }\quad u\in C_c^\infty\left(\overline{\R^{N+1}_+}\right). 
		\end{align*}
	\end{itemize} 
\end{lem}
{\sc{Proof.}}  The proof  follows   by using  Lemma \ref{Mazya 2} with  $\beta=\frac {c} q$ and, for claim (i),  by imposing  also $\frac 1 q=\frac 1 {2^\ast}$. 
\qed

	\begin{prop}\label{Sobolev2b}
		Let    $-1<c<0$ and $r>0$. The following properties hold.
		\begin{itemize}
			\item[(i)]If $N>1$ then  one has for some positive constant $C=C(r)$
			\begin{align*}
				\|u\|_{L^{2^\ast}_{c}(\R^{N+1}_+)}\leq C\, \|\nabla u\|_{L^2_{c}(\R^{N+1}_+)},\qquad \text{for all }\quad u\in H^1_c\left(\R^{N+1}_+\right),\quad \mbox{supp}\,u\subseteq \R^N\times [0,r].
			\end{align*}
			\item[(ii)]If $N=1$ then for every $q\geq 2$ there exists a  positive constant $C=C(r)$ such that  \begin{align*}
				\|u\|_{L^{q}_{c}(\R^{N+1}_+)}\leq C\, \|\nabla u\|_{L^2_{c}(\R^{N+1}_+)},\qquad \text{for all }\quad u\in H^1_c\left({\R^{N+1}_+}\right),\quad \mbox{supp}\,u\subseteq \R^N\times [0,r].
			\end{align*}
		\end{itemize} 
	\end{prop}
	{\sc{Proof.}} For simplicity we can assume $r=1$. Using the  density of  $\mathcal C$,  we can assume without any loss of generality $u\in\mathcal C$.  The proof then follows by (i) and (ii) of Lemma \ref{Sobolev2} after noticing that if $y\leq 1$ then  in (i) one has   $y^{-\frac c{N+1}}\leq 1$ and in (ii), $ y^{\frac{c+2}q-\frac c 2}\leq 1$.
	\qed 
	
	We can prove some  local  Gagliardo-Nirenberg inequalities.
	\begin{prop}\label{Sobolev3b}
		Let    $-1<c<0$, $r>0$  and $q\in (2,\infty)$ such that $\frac{1}{q}\geq \frac{1}{2}-\frac{1}{N+1}$.  Then for any $u\in H^1_c\left(\R^{N+1}_+\right)$ such that  $\mbox{supp}\,u\subseteq \R^N\times [0,r]$ one has  
		\begin{align*}
			\|u\|_{L^q_c(\R^{N+1}_+)}\leq C\, \|\nabla u\|_{L^2_{c}(\R^{N+1}_+)}^{(N+1)\left(\frac{1}{2}-\frac{1}{q}\right)}\| u\|_{L^2_c(\R^{N+1}_+)}^{1-(N+1)\left(\frac{1}{2}-\frac{1}{q}\right)}.
		\end{align*} 
	\end{prop}
	{\sc{Proof.}}  For simplicity we can assume $r=1$. The  claim for $N>1$ follows from (i) of Proposition \ref{Sobolev2b} using  Holder's inequality as in Proposition \ref{CKN2}.
	
	Let us prove the claim for $N=1$. Using  Lemma \ref{Mazya 1} with  $q=2$, $\beta=\alpha=\frac {c} 2$, we obtain 
	\begin{align*}
		\| u\|_{L^2_c(\R^{N+1}_+)}\leq C\|y^{-\frac c 2}\nabla u\|_{L^1_c(\R^{N+1}_+)},\quad \forall u\in C_c^\infty\left(\overline{\R^{N+1}_+}\right).
	\end{align*}
	In particular since $c<0$, we get
	\begin{align*}
		\| u\|_{L^2_c(\R^{N+1}_+)}\leq C\|\nabla u\|_{L^1_c(\R^{N+1}_+)},\quad \forall u\in C_c^\infty\left(\overline{\R^{N+1}_+}\right),\quad \mbox{supp}\,u\subseteq \R^N\times [0,1].
	\end{align*}
	Let $m\in\N$; applying the latter inequality to $|u|^{m-1}u$ instead of $u$ and using Holder's inequality  we obtain
	\begin{align*}
		\| u\|^m_{L^{2m}_c(\R^{N+1}_+)}&\leq Cm\|\nabla u|u|^{m-1}\|_{L^1_c(\R^{N+1}_+)}\leq Cm\|\nabla u\|_{L^2_c(\R^{N+1}_+)}\|u\|^{m-1}_{L^{2m-2}_c(\R^{N+1}_+)}
	\end{align*}
	that is 
	\begin{align}\label{Induzione1}
		\| u\|_{L^{2m}_c(\R^{N+1}_+)}& \leq (Cm)^{\frac 1 m}\|\nabla u\|_{L^2_c(\R^{N+1}_+)}^{\frac 1 m}\|u\|^{1-\frac 1 m}_{L^{2(m-1)}_c(\R^{N+1}_+)}
	\end{align}
	which is valid for any $u\in C_c^\infty\left(\overline{\R^{N+1}_+}\right)$ with $ \mbox{supp}\,u\subseteq \R^N\times [0,1]$. The inequality above and  an easy inductive argument then proves that, for some positive  constant $C=C(m)>0$, one has 
	\begin{align}\label{Induzione2}
		\| u\|_{L^{2m}_c(\R^{N+1}_+)}& \leq C \|\nabla u\|_{L^2_c(\R^{N+1}_+)}^{1-\frac 1 m}\|u\|^{\frac 1 m}_{L^{2}_c(\R^{N+1}_+)},\quad \forall u\in C_c^\infty\left(\overline{\R^{N+1}_+}\right),\quad \mbox{supp}\,u\subseteq \R^N\times [0,1]
	\end{align}
	which is the required claim for $q=2m$. Indeed for $m=1$ the latter inequality is simply $	\| u\|_{L^{2}_c(\R^{N+1}_+)} \leq C \|u\|_{L^{2}_c(\R^{N+1}_+)}$ which is obviously true. Let us suppose \eqref{Induzione2} true for $m=k$; then applying  \eqref{Induzione1} we get  for any $u\in C_c^\infty\left(\overline{\R^{N+1}_+}\right)$ with $ \mbox{supp}\,u\subseteq \R^N\times [0,1]$ and for some positive constant $C$ depending on $k$
\begin{align*}
		\| u\|_{L^{2(k+1)}_c(\R^{N+1}_+)}& \leq C\|\nabla u\|_{L^2_c(\R^{N+1}_+)}^{\frac 1 {k+1}}\|u\|^{1-\frac 1 {k+1}}_{L^{2k}_c(\R^{N+1}_+)}\\[1ex]
		& \leq C\|\nabla u\|_{L^2_c(\R^{N+1}_+)}^{\frac 1 {k+1}}\left(C \|\nabla u\|_{L^2_c(\R^{N+1}_+)}^{1-\frac 1 k}\|u\|^{\frac 1 k}_{L^{2}_c(\R^{N+1}_+)}\right)^{1-\frac 1 {k+1}}\\[1ex]
		&= C\|\nabla u\|_{L^2_c(\R^{N+1}_+)}^{1-\frac 1 {k+1}}\|u\|^{\frac 1 {k+1}}_{L^{2}_c(\R^{N+1}_+)}
\end{align*}
	which is the claim for $m=k+1$.
	
	 The claim for a general $q>2$  then  follows by interpolating $L^q_c(\R^{N+1}_+)$ between $L^{2m}_c(\R^{N+1}_+)$ and $L^{2m+2}_c(\R^{N+1}_+)$   for a suitable $m\in \N$.
	\qed 
	

\section{Kernel estimates}\label{Section Kernel}
In this section we deduce kernel estimates for the  degenerate   operator 
\begin{equation*} 
	\mathcal L=\Delta_x u+2a\cdot \nabla_xD_yu+ B_yu 
\end{equation*}
of Section \ref{L2}. We use the standard machinery of \cite[Chapter 6]{Ouhabaz},  relying on Gagliardo-Nirenberg type inequalities for deducing the ultracontractivity of the semigroup, which gives bounds on the diagonal, and on  Davies' method to add the gaussian term. Some difficulties appear  since the semigroup $e^{t\mathcal L}$ is not ultracontractive when $c<0$: this is due to the fact that the Sobolev immersions of  Proposition \ref{CKN2} fail to hold, in the whole space, for $c<0$ due to the singularity of the operator at the boundary  $y=0$. To overcome this obstruction, we introduce then a further  weight $\phi$ and an equivalent operator $\tilde{\mathcal L}$ which allow us to take into account the degeneracy of the operator at the boundary. This introduction is essential for $c<0$ since allows to prove the  Gagliardo-Nirenberg inequalities of Proposition \ref{Prop G-N ineq nu} which hold in the whole half-space and for any $c>-1$; when $c\geq 0$, this approach, although not mandatory,  is needed in order to deduce better heat kernel estimates (see Remark \ref{Oss c>0}).

\begin{defi}\label{Definizione phi}
	Let $\eta\in C^\infty[0,\infty)$ satisfying $0\leq \eta\leq 1$, $\eta=1$ on $[0,1/2]$ 
	and $\eta=0$ on $[2,\infty)$ and  let us 
	define the function $\phi\in C^\infty[0,+\infty)$ by
	\begin{equation}\label{def.phi}
		\phi(y):=\eta(y)+(1-\eta(y))y^{-\frac{c}{2}}.
	\end{equation}
\end{defi}
 By elementary computations, $\phi$ satisfies for some positive constant $C_0>0$, see also \cite[Lemma 3.3]{met-calv-negro-spina}:
\begin{align}\label{Prop-phi}
	\phi(y)\simeq y^{-\frac{c}{2}} (1\wedge y)^{\frac{c}{2}},\qquad 
	\left|D_y\phi\right|\leq C_0y ^{-1}\chi_{\{y\geq 1\}}\,\phi,\qquad |B\phi|\leq C_0 y ^{-2}\chi_{\{y\geq 1\}} \,\phi.
\end{align}
In particular 	$\left|D_y\phi\right|\leq C_0\,\phi$, $|B\phi|\leq C_0 \,\phi$.
\medskip 

\begin{defi}\label{defi tilde L}
 We introduce the Hilbert space $L^2_{\nu}:=L^2(\R^{N+1}_+,d\nu)$, 
and  consider the isometry 
\begin{align} \label{deffi}
	\Phi :L^2(\R^{N+1}_+,d\nu)\to L^{2}_c\left(\R^{N+1}_+\right), \qquad  f\mapsto \phi\,f
\end{align}
where by \eqref{Prop-phi}
\begin{align}\label{est measure}
	d\nu=\phi^2\,y^c dx\ dy\simeq (y\wedge 1)^c\,dxdy.
\end{align}
The map $\Phi$ defines a similar operator $\tilde {\mathcal L}=\Phi^{-1} \mathcal L \Phi $ which acts on $L^2(\R^{N+1}_+,d\nu)$.
\end{defi}
 In what follows we prove that the semigroup generated by $\tilde{\mathcal L}$ is  ultracontractive. We need some preparation. Let us consider the  symmetric form on  $L^2_\nu$ induced by $\Phi$ and  defined by
\begin{align*}
\mathfrak{b}(u,v)
:=\mathfrak{a}(\phi u,\phi v) \qquad\quad 
D(\mathfrak{b})=\Phi^{-1}\left(D(\mathfrak{a})\right)=\frac 1 \phi H^1_c.
\end{align*}
$\mathfrak{b}$ inherits, by construction, all the properties of $\mathfrak{a}$ stated in Propositions \ref{prop-form},  \ref{generation L2}. In particular, since $\mathcal C$ is preserved by $\Phi$, $\mathcal C$ is dense in $D(\mathfrak{b})$. Moreover $\mathfrak{b}$ is formally associated to $-\tilde{\mathcal L}$.  By \eqref{Prop-phi}, one  can recognize  that $D(\mathfrak{b})$ coincides with the Hilbert space 
 $$H^{1}_{\nu}\stackrel{def}{=}\overline{ \mathcal C}^{\|\cdot\|_{H^1_\nu}},\qquad 	\left\langle u, v\right\rangle_{H^1_{\nu}}:= \left\langle u, v\right\rangle_{L^2_{\nu}}+\left\langle \nabla u, \nabla v\right\rangle_{L^2_\nu}$$
where the closure is taken with respect to the norm
\begin{align*}
\|u\|_{H^1_{\nu}}:=\sqrt{ \|\nabla u\|^2_{L^2_{\nu}}+\|u\|^2_{L^2_\nu}}.
\end{align*}
We clarify the latter result in the following propositions.

\begin{lem}\label{Equivalent b}
	For any $u,v \in D(\mathfrak{b})$ one has  
	\begin{align} 
		\mathfrak{b}(u,v)
		&= \left\langle \nabla u, \nabla v\right\rangle_{L^2_\nu}+2\left\langle D_y u,a\cdot \nabla_x v\right\rangle_{L^2_\nu}  -\left\langle \frac{B\phi}{\phi}\,u,  v\right\rangle_{L^2_{\nu}}+2\left\langle \frac{D_y\phi}{\phi}\,u, a\cdot \nabla_x v\right\rangle_{L^2_{\nu}}
	\end{align}
\end{lem}
{\sc{Proof.}}
By density we can assume $u,v\in\mathcal C$; then by the definition of $\mathfrak{a}$ we have 
\begin{align*}
	 \mathfrak{a}(\phi u,\phi v)=&	\int_{\R^{N+1}_+} \langle \nabla (\phi u), \nabla (\phi \overline{v})\rangle\,y^{c} dx\,dy+2\int_{\R^{N+1}_+} D_y(\phi u)\, a\cdot\nabla_x\overline{v}\phi \,y^{c} dx\,dy\\[1ex]
	=&\int_{\R^{N+1}_+} \langle \nabla u, \nabla  \overline{v}\rangle\,\phi^2y^{c} dx\,dy+2\int_{\R^{N+1}_+} D_yu\, a\cdot\nabla_x\overline{v}\phi^2 \,y^{c} dx\,dy\\[1ex]
	&+	\int_{\R^{N+1}_+} \Big[
	D_y \phi D_y u\overline v\phi 
	+
	D_y \phi D_y\overline vu\phi
	+
	D_y \phi D_y\phi u\overline v\Big]\,y^{c} dx\,dy\\[1ex]
	&+2\int_{\R^{N+1}_+} D_y \phi u\, a\cdot\nabla_x\overline{v}\phi \,y^{c} dx\,dy
\end{align*}
which is equivalent to
\begin{align*}
	\mathfrak{a}(\phi u,\phi v)
	&= \left\langle \nabla u, \nabla v\right\rangle_{L^2_\nu}+2\left\langle D_y u,a\cdot \nabla_x v\right\rangle_{L^2_\nu}\\[1ex]
	&+	\int_{\R^{N+1}_+}D_y \phi D_y\left(u\overline v \phi\right) \,y^{c} dx\,dy+2\left\langle \frac{D_y\phi}{\phi}\,u, a\cdot \nabla_x v\right\rangle_{L^2_{\nu}}.
\end{align*}
The claim then follows after observing that using \eqref{Prop-phi} and integrating by parts we get  
\begin{align*}
		\int_{\R^{N+1}_+}D_y \phi D_y\left(u\overline v \phi\right) \,y^{c} dx\,dy&=-		\int_{\R^{N+1}_+}(B \phi) u\overline v \,\phi \,y^{c} dx\,dy+\Big[D_y\phi\, \phi\, u\overline v\,y^c\Big]_0^{+\infty}\\[1ex]
		&=-		\int_{\R^{N+1}_+}(B \phi) u\overline v\, \phi \,y^{c} dx\,dy.
\end{align*}
\qed

\begin{prop}
The following two norms are equivalent on $\mathcal C$:
	\begin{align*}
		\|u\|_{\mathfrak{b}}=\|\phi u\|_{\mathfrak{a}}=\sqrt{\Rp \mathfrak{a}(\phi u,\phi u)+\|\phi u\|_{L^2_c}},\qquad 
		\|u\|_{H^1_{\nu}}= \sqrt{\|\nabla u\|^2_{L^2_{\nu}}+\|u\|^2_{L^2_\nu}}
	\end{align*}
In particular   $D(\mathfrak{b})=H^1_\nu$.
\end{prop}
{\sc{Proof.}}
  We start by proving  that for some positive constant $C$ one has $\|u\|_{\mathfrak{b}}\leq C 
\|u\|_{H^1_{\nu}}$ for any $u\in\mathcal C$. Since $\|\phi u\|_{L^2_c}=\| u\|_{L^2_\nu}$, it is enough to prove that
\begin{align*}
\Rp \mathfrak{b}(u,u)=	\Rp \mathfrak{a}(\phi u,\phi u)\leq C 
	\|u\|_{H^1_{\nu}}^2,\qquad u\in\mathcal C.
\end{align*}
Let  $C_0$ be the constant in \eqref{Prop-phi}. Then using  Lemma \ref{Equivalent b},  the property of $\phi$ in \eqref{Prop-phi} and Young's inequality, we have 
\begin{align*}
	\Rp \mathfrak{b}(u,u)&= \left\langle \nabla u, \nabla u\right\rangle_{L^2_\nu}+2\left\langle D_y u,a\cdot \nabla_x u\right\rangle_{L^2_\nu}  -\left\langle \frac{B\phi}{\phi}\,u,  u\right\rangle_{L^2_{\nu}}+2\left\langle \frac{D_y\phi}{\phi}\,u, a\cdot \nabla_x u\right\rangle_{L^2_{\nu}}\\[1.5ex]
	&\leq  \|\nabla u\|^2_{L^2_\nu}+|a|\|\nabla u\|^2_{L^2_\nu}  +C_0\|u\|^2_{L^2_{\nu}}+C_0 |a|\left(\|u\|^2_{L^2_{\nu}}+\|\nabla_x u\|^2_{L^2_{\nu}}\right)\\[1.5ex]
	&\leq \left(1+|a|(1+C_0)\right)\|\nabla u\|^2_{L^2_\nu}+C_0(1+|a|)\,\| u\|^2_{L^2_\nu}
\end{align*}
which is the required claim.\\
On the other hand, with $C_0$ as in \eqref{Prop-phi} and using  Lemma \ref{Equivalent b},  the property of $\phi$ in \eqref{Prop-phi} and Young's inequality again we have also, for every $\eps>0$,
\begin{align*}
	\Rp \mathfrak{b}(u,u)&= \left\langle \nabla u, \nabla u\right\rangle_{L^2_\nu}+2\left\langle D_y u,a\cdot \nabla_x u\right\rangle_{L^2_\nu}  -\left\langle \frac{B\phi}{\phi}\,u,  u\right\rangle_{L^2_{\nu}}+2\left\langle \frac{D_y\phi}{\phi}\,u, a\cdot \nabla_x u\right\rangle_{L^2_{\nu}}\\[1.5ex]
	&\geq  (1-|a|)\|\nabla u\|^2_{L^2_\nu}-C_0\|u\|^2_{L^2_{\nu}}-\frac{C_0}{\eps}|a|\|u\|^2_{L^2_{\nu}}-C_0|a|\eps\|\nabla_x u\|^2_{L^2_{\nu}}.\\[1.5ex]
	\end{align*}
	Choosing $\eps$ small enough such that $(1-|a|)-C_0|a|\epsilon>0$, it follows that
	$$\|\nabla u\|^2_{{L^2_\nu}}\leq C\left(\Rp\mathfrak{b}(u,u)+\| u\|^2_{{L^2_\nu}}\right)$$ for some positive constant $C$. \qed


By similarity and by Proposition \ref{generation L2}, 	$ \tilde{\mathcal L}$ generates a contractive analytic semigroup of angle $\frac{\pi}{2}-\arctan \frac{|a|}{1-|a|}$  in $L^2_{\nu}$.  However, in order to deduce pointwise kernel estimates by Davies' method, we need to introduce perturbed forms as follows.

\begin{defi}Let 
\begin{equation}\label{distortion}
\Psi:=
\Big\{\psi\in C^\infty(\overline{\R^{N+1}_+})\cap L^\infty(\R^{N+1}_+)\;;|\nabla \psi|\leq 1;\quad |D^2\psi|\leq 1\Big\}. 
\end{equation}
We introduce the quadratic forms $\{\mathfrak{b}_{\alpha\psi}\}_{\alpha\in\R,\psi\in\Psi}$ 
in $L^2_{\nu}$ defined by
\begin{align*}
\mathfrak{b}_{\alpha\psi}(u,v)
=\mathfrak{a}(e^{\alpha\psi}\phi u,e^{-\alpha\psi}\phi v), \qquad D(\mathfrak{b}_{\alpha\psi})=D(\mathfrak{b}).
\end{align*}
\end{defi}

\begin{os}\label{equiv distance}
	By standard computation one can recognize that 
	\begin{align*}
		d(z_1,z_2)\stackrel{def}=\sup_{\psi\in\Psi} \{\psi(z_1)-\psi(z_2)\},\qquad \forall z_1,z_2\in\R^{N+1}_+
	\end{align*}
defines a distance in $\R^{N+1}$ which is actually equivalent to the euclidean one i.e.
	\begin{align*}
		|z_1-z_2|\simeq\sup_{\psi\in\Psi} \{\psi(z_1)-\psi(z_2)\},\qquad \forall z_1,z_2\in\R^{N+1}_+.
	\end{align*}
We refer the reader to \cite[pages 200-203]{Robinson-elliptic} or to \cite[Corollary 3.2.8]{davies} for the proof of this classic result.
\end{os}

Since $\Phi$ is an isometry from $L^2_{\nu}$ onto $L^2_{c}$  which  maps $D(\mathfrak{b})$ onto $D(\mathfrak{a})$ and since
$v\in L^2_{\nu}\to e^{\pm\alpha\psi}v\in L^2_{\nu}$ are bounded operators 
preserving $D(\mathfrak{b})$, the following result  easily follows.
\begin{prop}\label{basic}
	For every $\alpha\in\R$ and $\psi\in \Psi$ the operator  
	$$\mathcal B_{\alpha, \psi}u=\phi^{-1}e^{-\alpha\psi}\mathcal L(\phi e^{\alpha\psi}u)$$ 
	with domain $D(\mathcal B_{\alpha\psi})=\phi^{-1}e^{-\alpha\psi}(D(\mathcal L))$ is associated 
	to the  form $\mathfrak{b}_{\alpha\psi}$. 
	Moreover $\mathcal B_{\alpha, \psi}$   generates a bounded analytic semigroups $e^{z \mathcal B_{\alpha, \psi}},\;z\in\Sigma_{\frac{\pi}{2}-\arctan \frac{|a|}{1-|a|}}$  in $L^2_{\nu}$ which satisfies
	\begin{align}\label{sem-pert equiv}
		e^{z\mathcal B_{\alpha\psi}}f= \phi^{-1}e^{-\alpha\psi}e^{z\mathcal L}
		\left(\phi e^{\alpha\psi}f\right), \quad f \in L^2_\nu.
	\end{align}
\end{prop}

\begin{prop}\label{accretive}
	For every $\alpha\in\R$ and $\psi\in\Psi$, $\mathfrak{b}_{\alpha\psi}$ 
	is quasi-accretive in $L^2_{\nu}$  and for any $0<\epsilon<1$ there exists some positive constant $C_\epsilon>0$ dependent only on $\epsilon, a, c$ such that 
	\begin{align*}
		{\Rp \mathfrak{b}_{\alpha\psi}}(u,u)
		\geq& \,\epsilon(1-|a|)\|\nabla u\|^2_{L^2_\nu}-C_\epsilon\left(\alpha^2+1\right)\|u\|^2_{L^2_\nu}.
	\end{align*}
Moreover, $e^{t\mathcal B_{\alpha\psi}}$ is positive and one has 
\begin{align*}
			\|e^{t\mathcal B_{\alpha\psi}}f\|_{\infty}
		\leq 
		e^{C_\epsilon(1+\alpha^2)t}\|f\|_{\infty},
		\quad f\in L^2_{\nu}\cap L^\infty (\R^{N+1}_+).
\end{align*}

\end{prop}
{\sc Proof.}\  Since the form is real {(i.e. $\mathfrak{b}_{\alpha\psi}(u,v)\in \R$ if $u,\ v\in\R$)}, {one has $\Rp\mathfrak{b}_{\alpha\psi}(u,u)=\mathfrak{b}_{\alpha\psi}(\Rp u,\Rp u)+\mathfrak{b}_{\alpha\psi}(\Ip u,\Ip u)$} so we can restrict the analysis,  without any loss of generality, to real  $u\in D(\mathfrak{b})$.
Using Lemma \ref{Equivalent b} and a straightforward calculation we  get
\begin{align*}
	\mathfrak{b}_{\alpha\psi}(u,u)
	=&
	\left\langle \nabla (e^{\alpha\psi}u), \nabla (e^{-\alpha\psi}u)\right\rangle_{L^2_\nu}+2\left\langle D_y (e^{\alpha\psi}u),a\cdot \nabla_x (e^{-\alpha\psi}u)\right\rangle_{L^2_\nu} 
	 \\[1ex]
	&-\left\langle \frac{B\phi}{\phi}\,e^{\alpha\psi}u,  e^{-\alpha\psi}u\right\rangle_{L^2_{\nu}}+2\left\langle \frac{D_y\phi}{\phi}\,e^{\alpha\psi}u, a\cdot \nabla_x (e^{-\alpha\psi}u)\right\rangle_{L^2_{\nu}}
	\\[1ex]
	=&
	\left\|\nabla u\right\|^2_{L^2_\nu}-\alpha^2	\left\langle \nabla \psi u,\nabla \psi u\right\rangle_{L^2_\nu}-2\alpha^2	\left\langle D_y \psi u,a\cdot\nabla_x \psi u\right\rangle_{L^2_\nu}+2\left\langle D_y u,a\cdot\nabla_x  
	u\right\rangle_{L^2_\nu}
	\\[1ex]
	&-2\alpha	\left\langle D_y  u,a\cdot\nabla_x\psi  
	u\right\rangle_{L^2_\nu}+2\alpha	\left\langle D_y \psi u,a\cdot\nabla_x 
	u\right\rangle_{L^2_\nu}
	\\[1ex]&
	-\left\langle \frac{B\phi}{\phi} u,u\right\rangle_{L^2_\nu}-2\alpha \left\langle \frac{D_y\phi}{\phi} u,a\cdot\nabla_x\psi u\right\rangle_{L^2_\nu}+2\left\langle \frac{D_y\phi}{\phi} u,a\cdot\nabla_x u\right\rangle_{L^2_\nu}
\end{align*}
Using the properties of $\phi$ in \eqref{Prop-phi} we obtain
\begin{align}\label{accretive1}
	\nonumber\mathfrak{b}_{\alpha\psi}(u,u)
	\geq &	\left\|\nabla u\right\|^2_{L^2_\nu}-\alpha^2\|u\|^2_{L^2_\nu}-2\alpha^2|a|\|u\|^2_{L^2_\nu}-2|a|\|D_yu\|_{L^2_\nu}\|\nabla_xu\|_{L^2_\nu}\\[1ex]\nonumber
	&-2|\alpha||a|\|u\|_{L^2_\nu}\|D_yu\|_{L^2_\nu}-2|\alpha||a|\|u\|_{L^2_\nu}\|\nabla_xu\|_{L^2_\nu}\\[1ex]
	&-C_0\|u\|^2_{L^2_\nu}-2|\alpha||a|C_0\|u\|^2_{L^2_\nu}-2C_0|a|\|u\|_{L^2_\nu}\|\nabla_xu\|_{L^2_\nu}.
\end{align}
Using the Cauchy-Schwarz inequality we get for any $\epsilon_1>0$
\begin{align*}
	2|\alpha|\|u\|_{L^2_\nu}\|D_yu\|_{L^2_\nu}=2\left(\frac{|\alpha|}{\sqrt{\epsilon_1}}\|u\|_{L^2_\nu}\right)\left(\sqrt{\epsilon_1}\|D_yu\|_{L^2_\nu}\right)\leq \frac{\alpha^2}{\epsilon_1}\|u\|^2_{L^2_\nu}+\epsilon_1\|D_yu\|^2_{L^2_\nu}
\end{align*}
and similarly $2|\alpha|\|u\|_{L^2_\nu}\|\nabla_xu\|_{L^2_\nu}\leq \frac{\alpha^2}{\epsilon_1}\|u\|^2_{L^2_\nu}+\epsilon_1\|\nabla_xu\|^2_{L^2_\nu}$ which implies
\begin{align*}
	2|\alpha|\|u\|_{L^2_\nu}\|D_yu\|_{L^2_\nu}+2|\alpha|\|u\|_{L^2_\nu}\|\nabla_xu\|_{L^2_\nu}\leq 2\frac{\alpha^2}{\epsilon_1}\|u\|^2_{L^2_\nu}+\epsilon_1\|\nabla u\|^2_{L^2_\nu}.
\end{align*}
Analogously we get $2|a|\|D_yu\|_{L^2_\nu}\|\nabla_xu\|_{L^2_\nu}\leq |a|\|\nabla u\|^2_{L^2_\nu}$ and  for any $\epsilon_2>0$
\begin{align*}
	2\|u\|_{L^2_\nu}\|\nabla_xu\|_{L^2_\nu}\leq\frac{1}{\epsilon_2}\|u\|^2_{L^2_\nu}+\epsilon_2\|\nabla_xu\|^2_{L^2_\nu}
\end{align*}
Inserting the latter estimates into \eqref{accretive1} yields
\begin{align*}
	\mathfrak{b}_{\alpha\psi}(u,u)
	\geq& \Big[1-|a|(1+\epsilon_1+C_0\epsilon_2)\Big]\|\nabla u\|^2_{L^2_\nu}\\[1ex]
	&-\Big[\alpha^2\left(1+2|a|+\frac 2 {\epsilon_1}|a|+\frac{C_0}{\epsilon_2}\right)+2|\alpha||a|C_0+C_0+\frac{C_0|a|}{\epsilon_2}\Big]\|u\|^2_{L^2_\nu}.
\end{align*}
Let us fix now $0<\epsilon<1$. Since $|a|<1$ we can choose $\epsilon_1,\epsilon_2$ sufficiently small such that $$1-|a|(1+\epsilon_1+C_0\epsilon_2)\geq \epsilon\left(1-|a|\right).$$
Consequently, using the elementary inequality  $\alpha^2+c_1|\alpha|+c_2 \le 2\alpha^2+c_3$ for  a suitable $c_3>0$, we get  for some positive constant $C_1>0$ dependent only on $\epsilon, a$
\begin{align*}
	\mathfrak{b}_{\alpha\psi}(u,u)
	\geq& \epsilon(1-|a|)\|\nabla u\|^2_{L^2_\nu}-C_1\left(\alpha^2+1\right)\|u\|^2_{L^2_\nu}.
\end{align*}

The proof of the remaining claim follows similarly as in the proof of  Proposition \ref{generation L2}. Indeed the positivity follows by \cite[Theorem 2.6]{Ouhabaz} after observing that, if $u\in H^1_{\nu}$, $u$ real, then $u^+\in H^1_{\nu}$ and, by Proposition \ref{generation L2},
\begin{align*}
	\mathfrak{b}_{\alpha\psi}(u^+,u^-)
	&=\mathfrak{a}(e^{\alpha\psi}\phi u^+,e^{-\alpha\psi}\phi u^-)=\mathfrak{a}((e^{\alpha\psi}\phi u)^+,(e^{-\alpha\psi}\phi u)^-)
=0.
\end{align*}
Finally, the quasi-$L^\infty$-contractivity follows by \cite[Corollary 2.17]{Ouhabaz}  after proving that if $0\leq u\in H^1_{\nu}$  then, observing that  $1\wedge u, (u-1)^+\in H^1_{\nu}$, one has  $$\mathfrak{b}_{\alpha\psi}(1\wedge u,(u-1)^+)+C_2 (1\wedge u,(u-1)^+)\geq 0$$ for some positive constant $C_2>0$ independent of $\alpha, \psi$. After this the required claim of the Proposition follows after taking $C_\epsilon=\max\{C_1,C_2\}$.  
Indeed, let $0\leq u\in H^1_{\nu}$ and, using the density of $\mathcal C$ in \eqref{defC},  let us assume, without any loss of generality, $u\in \mathcal C$. Then 

\begin{align*}
	\mathfrak{b}&_{\alpha\psi}(1\wedge u,(u-1)^+)\\[1ex]
	=&
	\left\langle \nabla (e^{\alpha\psi}(1\wedge u)), \nabla (e^{-\alpha\psi}(u-1)^+)\right\rangle_{L^2_\nu}+2\left\langle D_y (e^{\alpha\psi}(1\wedge u)),a\cdot \nabla_x (e^{-\alpha\psi}(u-1)^+)\right\rangle_{L^2_\nu} 
	\\[1ex]
	&-\left\langle \frac{B\phi}{\phi}\,e^{\alpha\psi}(1\wedge u),  e^{-\alpha\psi}(u-1)^+\right\rangle_{L^2_{\nu}}+2\left\langle \frac{D_y\phi}{\phi}\,e^{\alpha\psi}(1\wedge u), a\cdot \nabla_x (e^{-\alpha\psi}(u-1)^+)\right\rangle_{L^2_{\nu}}
	\\[1ex]
	=&
		\left\langle \nabla(1\wedge u),\nabla  (u-1)^+\right\rangle_{L^2_\nu}-\alpha^2	\left\langle \nabla \psi (1\wedge u),\nabla \psi (u-1)^+\right\rangle_{L^2_\nu}
		\\[1ex]
		&+\alpha \left\langle\nabla\psi (1\wedge u),\nabla(u-1)^+\right\rangle_{L^2_\nu}-\alpha \left\langle\nabla (1\wedge u),\nabla\psi (u-1)^+\right\rangle_{L^2_\nu}
		\\[1ex]
		&-2\alpha^2	\left\langle D_y \psi (1\wedge u),a\cdot\nabla_x \psi (u-1)^+\right\rangle_{L^2_\nu}+2\left\langle D_y (1\wedge u),a\cdot\nabla_x  
	(u-1)^+\right\rangle_{L^2_\nu}
	\\[1ex]
	&-2\alpha	\left\langle D_y  (1\wedge u),a\cdot\nabla_x\psi  
	(u-1)^+\right\rangle_{L^2_\nu}+2\alpha	\left\langle D_y \psi (1\wedge u),a\cdot\nabla_x 
	(1-u)^+\right\rangle_{L^2_\nu}
	\\[1ex]&
	-\left\langle \frac{B\phi}{\phi} (1\wedge u),(u-1)^+\right\rangle_{L^2_\nu}-2\alpha \left\langle \frac{D_y\phi}{\phi} (1\wedge u),a\cdot\nabla_x\psi u\right\rangle_{L^2_\nu}+2\left\langle \frac{D_y\phi}{\phi} (1\wedge u),a\cdot\nabla_x (1-u)^+\right\rangle_{L^2_\nu}.
\end{align*}
Since $\nabla (1\wedge u)=\chi_{\{u<1\}}\nabla u$ 
and 
$\nabla (u-1)^+=\chi_{\{u>1\}}\nabla u$, one has
\begin{align*}
	\mathfrak{b}_{\alpha\psi}(1\wedge u,(u-1)^+)
	= &-\alpha^2	\left\langle \nabla \psi ,\nabla \psi (u-1)^+\right\rangle_{L^2_\nu}+\alpha \left\langle\nabla\psi ,\nabla (u-1)^+\right\rangle_{L^2_\nu}\\[1ex]
	&-2\alpha^2	\left\langle D_y \psi ,a\cdot\nabla_x \psi (u-1)^+\right\rangle_{L^2_\nu}
	+2\alpha	\left\langle D_y \psi ,a\cdot\nabla_x 
	 (u-1)^+\right\rangle_{L^2_\nu}
	\\[1ex]&
	-\left\langle \frac{B\phi}{\phi} ,(u-1)^+\right\rangle_{L^2_\nu}-2\alpha \left\langle \frac{D_y\phi}{\phi} ,a\cdot\nabla_x\psi (u-1)^+\right\rangle_{L^2_\nu}+2\left\langle \frac{D_y\phi}{\phi},a\cdot\nabla_x (u-1)^+\right\rangle_{L^2_\nu}.
\end{align*}

 Using the properties of $\psi$ in \eqref{distortion} and  of $\phi$ in \eqref{Prop-phi} we get
 \begin{align*}
 	\mathfrak{b}_{\alpha\psi}(1\wedge u,(u-1)^+)
 	\geq  &\left(-\alpha^2-2|a|\alpha^2-C_0-2|\alpha|C_0|a|\right) \int_{\R^{N+1}_+} (u-1)^+d\nu\\[1ex]
 		&+\alpha \left\langle\nabla\psi ,\nabla (u-1)^+\right\rangle_{L^2_\nu}	+2\alpha	\left\langle D_y \psi ,a\cdot\nabla_x 
 	(u-1)^+\right\rangle_{L^2_\nu}+2\left\langle \frac{D_y\phi}{\phi},a\cdot\nabla_x (u-1)^+\right\rangle_{L^2_\nu}.
 \end{align*}
 
To estimate the second term  we  integrate by parts and we use again the properties of $\psi$ and  of $\phi$ (which implies in particular the estimate $\frac{|D_y\left(\phi^2 y^c\right)|}{\phi^2y^c}\leq C $, for some $C>0$) obtaining
\begin{align*}
\left|\alpha \left\langle\nabla\psi ,\nabla (u-1)^+\right\rangle_{L^2_\nu}\right|&=\left|\alpha \int_{\R^{N+1}_+}\nabla\psi \cdot\nabla (u-1)^+d\nu\right|\\[1ex]
&=\left|-\alpha \int_{\R^{N+1}_+}\left(\Delta\psi+\frac{D_y\psi D_y (\phi^2 y^c)}{\phi^2 y^c} \right) (u-1)^+d\nu\right|\\[1ex]
&\leq |\alpha|(N+1+C) \int_{\R^{N+1}_+} (u-1)^+d\nu.
\end{align*}
Similarly we get
\begin{align*}
\left|2\alpha	\left\langle D_y \psi ,a\cdot\nabla_x 
(u-1)^+\right\rangle_{L^2_\nu}\right|
	&\leq 2|\alpha||a| \int_{\R^{N+1}_+} (u-1)^+d\nu\\[1ex]
	2\left\langle \frac{D_y\phi}{\phi},a\cdot\nabla_x (u-1)^+\right\rangle_{L^2_\nu}
	&=0.
\end{align*}
Putting together all the above estimates and using again the elementary inequality  $\alpha^2+c_1|\alpha|+c_2 \le 2\alpha^2+c_3$, we finally get, for some positive constant $C>0$,
\begin{align*}
	\mathfrak{b}_{\alpha\psi}(1\wedge u,(u-1)^+)
& \geq -C\left(\alpha^2	+1 \right)\int_{\R^{N+1}_+} (u-1)^+d\nu.
\end{align*}

 \qed

By analogous computations we can get similar results for the adjoint form $\{\mathfrak{b}^\ast_{\alpha\psi}\}_{\alpha\in\R,\psi\in\Psi}$ 
in $L^2_{\nu}$ defined by
\begin{align*}
\mathfrak{b}^\ast_{\alpha\psi}(u,v)=\overline{\mathfrak{b}^\ast_{\alpha\psi}(v,u)}
=\mathfrak{a^*}(e^{-\alpha\psi}\phi u,e^{\alpha\psi}\phi v), \qquad D(\mathfrak{b^*}_{\alpha\psi})=D(\mathfrak{b})
\end{align*}
whose associated operator is the  adjoint operator of $\mathcal B_{\alpha, \psi}$
 in $L^2_\nu$ which we denote by $\mathcal B^*_{\alpha, \psi}$.

\begin{prop}\label{basic*}
	For every $\alpha\in\R$ and $\psi\in \Psi$ the operator  
	$$\mathcal B^*_{\alpha, \psi}u=\phi^{-1}e^{\alpha\psi}\mathcal L^*(\phi e^{-\alpha\psi}u)$$ 
	with domain $D(\mathcal B^*_{\alpha\psi})=\phi^{-1}e^{\alpha\psi}(D(\mathcal L^*))$ is associated 
	to the  form $\mathfrak{b}^\ast_{\alpha\psi}$. 
	Moreover $\mathcal B^*_{\alpha, \psi}$   generates a bounded analytic semigroups $e^{z \mathcal B^*_{\alpha, \psi}},\;z\in\Sigma_{\frac{\pi}{2}-\arctan \frac{|a|}{1-|a|}}$  in $L^2_{\nu}$ which satisfies
	\begin{align}\label{sem-pert equiv}
		e^{z\mathcal B^*_{\alpha\psi}}f= \phi^{-1}e^{\alpha\psi}e^{z\mathcal L^*}
		\left(\phi e^{-\alpha\psi}f\right), \quad f \in L^2_\nu.
	\end{align}
\end{prop}

As before, we can also prove that
\begin{prop}\label{accretive*}
	For every $\alpha\in\R$ and $\psi\in\Psi$, $\mathfrak{b^*}_{\alpha\psi}$ 
	is quasi-accretive in $L^2_{\nu}$  and for any $0<\epsilon<1$ there exist some positive constant $C_\epsilon>0$ dependent only on $\epsilon, a$ such that 
	\begin{align*}
		{\Rp\mathfrak{b}}_{\alpha\psi}^*(u,u)
		\geq& \,\epsilon(1-|a|)\|\nabla u\|^2_{L^2_\nu}-C_\epsilon\left(\alpha^2+1\right)\|u\|^2_{L^2_\nu}.
	\end{align*}
Moreover, $e^{t\mathcal B^*_{\alpha\psi}}$ is positive and one has 
				\begin{align*}
					\|e^{t\mathcal B^*_{\alpha\psi}}f\|_{\infty}
		\leq 
		e^{C_\epsilon(1+\alpha^2)t}\|f\|_{\infty},
		\quad f\in L^2_{\nu}\cap L^\infty (\R^{N+1}_+).
				\end{align*}		
\end{prop}

\subsection{Gagliardo-Nirenberg type inequalities for $\mathfrak{b}_{\alpha\psi}$ and $\mathfrak{b}^\ast_{\alpha\psi}$}

In this section we use the results of  Section \ref{Section weighted sobolev} to derive some  Gagliardo-Nirenberg inequalities for the space $H^1_\nu$ which are the main technical tools to prove upper bounds for the heat kernel of $\mathcal L$.

\begin{prop}\label{Prop G-N ineq nu}
	Assume  $c+1>0$. 
	Then for every $q\in(2,\infty)$ such that 
	$\frac{1}{q}\geq\frac{1}{2}-\frac{1}{N+1+c^+}$ one has  
	$$
	\|u\|_{L^q_{\nu}}
	\leq 
	C_q
	\|u\|_{H^1_\nu}^\theta
	\|u\|_{\nu}^{1-\theta},\qquad \text{for all}\quad u\in H^1_{\nu}.
	$$
	where  $C_q$ is a positive constant independent of $u$ and $\theta=\left (N+1+c^+\right )\left(\frac{1}{2}-\frac{1}{q}\right ) \in (0,1]$ ($\theta<1$ if $N=1$ and $c\leq 0$).
\end{prop}
{\sc Proof.}\ 
Let $\eta\in C^\infty[0,\infty)$ satisfying $\eta=1$ on $[0,1/2]$ 
and $\eta=0$ on $[2,\infty)$ and $\phi$ as defined in \eqref{def.phi}. 
We note preliminarily that by \eqref{est measure} we have
\begin{align*}
	d\nu\simeq  y^cdxdy,\quad \text{if \;} y\in {\rm supp}\,\eta,\qquad 	 d\nu\simeq dxdy,\quad \text{if \;} y\in {\rm supp}\,(1-\eta).
\end{align*}
We treat separately the case $c<0$ and $c\geq 0$. 

We start by the case $c<0$. From Lemmas \ref{Sobolev2b}  and \ref{GN}, respectively, we deduce,  for some positive constant $C=C_q>0$ which may depend on $q$ and which may vary in each occurrence, that 
\begin{align*}
	\|\eta u\|_{L^q_{\nu}}&\leq C\|\eta u\|_{L^q_{c}}\leq C_q \|\nabla(\eta  u)\|_{L^2_{c}}^{(N+1)\left(\frac{1}{2}-\frac{1}{q}\right)}
	\|\eta u\|_{L^2_{c}}^{1-(N+1)\left(\frac{1}{2}-\frac{1}{q}\right)}\\[1ex]
	&\leq 
	C_q\|\eta\|^{(N+1)\left(\frac{1}{2}-\frac{1}{q}\right)}_{W^{1,\infty}} \Big(\|u\|_{L^2_{\nu}}+\|\nabla u\|_{L^2_{\nu}}\Big)^{(N+1)\left(\frac{1}{2}-\frac{1}{q}\right)}
	\|u\|_{L^2_{\nu}}^{1-(N+1)\left(\frac{1}{2}-\frac{1}{q}\right)}
\end{align*}
and 
\begin{align*}
	\|(1-\eta)u\|_{L^q_{\nu}}
	&\leq 
	C\|(1-\eta)u\|_{q}
	\\
	&\leq 
	{C}_q\left\|\nabla((1-\eta)u)\right\|_{2}^{(N+1)\left(\frac{1}{2}-\frac{1}{q}\right)}
	\|(1-\eta)u\|_2^{1-(N+1)\left(\frac{1}{2}-\frac{1}{q}\right)}
	\\
	&\leq 
	{C}_q\|\eta\|^{(N+1)\left(\frac{1}{2}-\frac{1}{q}\right)}_{W^{1,\infty}} \Big(\|u\|_{L^2_{\nu}}+\|\nabla u\|_{L^2_{d\nu}}\Big)^{(N+1)\left(\frac{1}{2}-\frac{1}{q}\right)}
	\|u\|_{L^2_{\nu}}^{1-(N+1)\left(\frac{1}{2}-\frac{1}{q}\right)}.
\end{align*}
Combining the above estimates, we deduce 
\[
\|u\|_{L^q_{\nu}}\leq C
\Big(\|u\|_{L^2_{\nu}}+\|\nabla u\|_{L^2_{\nu}}\Big)^{(N+1)\left(\frac{1}{2}-\frac{1}{q}\right)}
\|u\|_{L^2_{\nu}}^{1-(N+1)\left(\frac{1}{2}-\frac{1}{q}\right)}
\]
which is the required claim.

We treat now, similarly,  the case $c\geq 0$.
 Observing that $\frac{1}{q}> \frac{1}{2}-\frac{1}{N+1+c}> \frac{1}{2}-\frac{1}{N+1}$, 
and applying Lemmas \ref{Sobolev2b}  and \ref{GN}, respectively, we get proceeding as before that 
\begin{align*}
	\|\eta u\|_{L^q_{\nu}}
	&\leq 
	C_q \Big(\|u\|_{L^2_{\nu}}+\|\nabla u\|_{L^2_{\nu}}\Big)^{(N+1+c)\left(\frac{1}{2}-\frac{1}{q}\right)}
	\|u\|_{L^2_{\nu}}^{1-(N+1+c)\left(\frac{1}{2}-\frac{1}{q}\right)}
\end{align*}
and
\begin{align*}
	\|(1-\eta)u\|_{L^q_{\nu}}
	&\leq 
	C_q \Big(\|u\|_{L^2_{\nu}}+\|\nabla u\|_{L^2_{\nu}}\Big)^{(N+1)\left(\frac{1}{2}-\frac{1}{q}\right)}
	\|u\|_{L^2_{\nu}}^{1-(N+1)\left(\frac{1}{2}-\frac{1}{q}\right)}
	\\[1ex]
	&\leq 
	C_q
	\Big(\|u\|_{L^2_{\nu}}+\|\nabla u\|_{L^2_{\nu}}\Big)^{(N+1+c)\left(\frac{1}{2}-\frac{1}{q}\right)}
	\|u\|_{L^2_{\nu}}^{1-(N+1+c)\left(\frac{1}{2}-\frac{1}{q}\right)}
\end{align*}
where in the last inequality we used the fact that  $(A+B)^{r}A^{1-r}$ is increasing in $r \in [0,1]$, for $A,B \ge 0$.
Combining the above estimates, we deduce \begin{align*}
\|u\|_{L^q_{\nu}}\leq C_q
\Big(\|u\|_{L^2_{\nu}}+\|\nabla u\|_{L^2_{\nu}}\Big)^{(N+1+c)\left(\frac{1}{2}-\frac{1}{q}\right)}
\|u\|_{L^2_{\nu}}^{1-(N+1+c)\left(\frac{1}{2}-\frac{1}{q}\right)}
\end{align*}
which is the required claim.
\qed

Using Propositions \ref{accretive} and \ref{accretive*} we can immediately reformulate the latter result in terms of $\mathfrak{b}_{\alpha\psi}$ and $\mathfrak{b}^\ast_{\alpha\psi}$.

\begin{cor}\label{Gagl-Nir-bpert}
		Assume  $c+1>0$. 
	Then for every $q\in(2,\infty)$ such that 
	$\frac{1}{q}\geq\frac{1}{2}-\frac{1}{N+1+c^+}$ one has for all  $u\in H^1_{\nu}=D(\mathfrak{b}_{\alpha\psi})=D(\mathfrak{b}^\ast_{\alpha\psi})$
	$$
	\|u\|_{L^q_{\nu}}
	\leq 
	C\Big(
	(1+\alpha^2)\|u\|_{\nu}^2
	+{\Rp \mathfrak{b}_{\alpha\psi}}(u,u)
	\Big)^{\frac{\theta}{2}}
	\|u\|_{\nu}^{1-\theta}
	$$
	and
	$$
	\|u\|_{L^q_{\nu}}
	\leq 
	C\Big(
	(1+\alpha^2)\|u\|_{\nu}^2
	+{\Rp\mathfrak{b}}_{\alpha\psi}^*(u,u)
	\Big)^{\frac{\theta}{2}}
	\|u\|_{\nu}^{1-\theta}
	$$
	where   $\theta=\left (N+1+c^+\right )\left(\frac{1}{2}-\frac{1}{q}\right ) \in (0,1]$ ($\theta<1$ if $N=1$ and $c\leq 0$) and $C>0$ is a positive constant which may depend only on $q,a,c$.
\end{cor}

\subsection{Derivation of kernel estimates for $e^{z{\mathcal L}}$ and  $e^{z{\mathcal L}^\ast}$}
We obtain kernel estimates  for $t>0$ via Davies' method. For  $x\in\R^N$, $y\in \R_+$,  we write  $z=(x,y)$ to denote by a generic point in $\R^{N+1}_+$. We consider preliminarily the operator $\tilde {\mathcal L}$  on $L^2(\R^{N+1}_+,d\nu)$ of 
Definition \ref{defi tilde L}; we also consider the perturbed operators $\mathcal B_{\alpha, \psi}$ of Proposition \ref{basic} and,  for clarity sake,  we recall that  
\begin{align}\label{relation tilde L pert}
\tilde {\mathcal L}=\mathcal B_{0, \psi},\qquad \mathcal B_{\alpha, \psi}u=e^{-\alpha\psi}\tilde {\mathcal L}(e^{\alpha\psi}u),\qquad e^{t\mathcal B_{\alpha, \psi}}f=e^{-\alpha\psi}e^{t\tilde{\mathcal L}}(e^{\alpha\psi}f).
\end{align}
\begin{lem} \label{realheat tilde}
	Let $c+1>0$, $|a|<1$. The semigroup $(e^{t\tilde{\mathcal L}})_{t> 0}$ consists of integral operators i.e. there exists $p_{\tilde{\mathcal L}}(t,\cdot,\cdot)\in L^\infty(\R^{N+1}_+\times\R^{N+1}_+)$ such that  for $t>0$, $z_1,\ z_2\in\R^{N+1}_+$
\begin{align*}
	e^{t\tilde{\mathcal  L}}f(z_1)=
\int_{\R^{N+1}_+}p_{\tilde{\mathcal L}}(t,z_1,z_2)f(z_2)\,d\nu,\quad f\in L^2_{\nu}.
\end{align*}
  Its heat kernel $p_{\tilde{\mathcal L}}$, written  with respect the measure $d\nu$, satisfies for some  $C, \kappa>0$, 
\begin{align}\label{Stima ptilde reale}
	0\leq p_{\tilde{\mathcal  L}}(t,z_1,z_2)
\leq C t^{-\frac{N+1+c^+}{2}}\exp\left(C t-\dfrac{|z_1-z_2|^2}{kt}\right).
\end{align}
\end{lem}

{\sc Proof.}\ We use Propositions \ref{accretive} and \ref{accretive*} with $\epsilon=\frac 1 2 $ and correspondingly we consider the forms 
$$\mathfrak{b}_{\alpha,\psi}+C(1+\alpha^2)\langle\, \, , \,\, \rangle_{L^2_\nu},\quad \mathfrak{b}^*_{\alpha,\psi}+C(1+\alpha^2)\langle\, \, , \,\, \rangle_{L^2_\nu}$$
which are  accretive in $L^2_\nu$  and $L^\infty$-contractive. 
  By Corollary \ref{Gagl-Nir-bpert} we can apply \cite[Theorem 6.2]{Ouhabaz} to get that \begin{align*}
  	\|e^{t\mathcal B_{\alpha,\psi}}\|_{L^2_{\nu}\to L^\infty}\leq Ct^{-\frac{N+1+c^+}{4}}e^{C(1+\alpha^2)t},\quad \|e^{t\mathcal B^*_{\alpha,\psi}}\|_{L^2_{\nu}\to L^\infty}\leq Ct^{-\frac{N+1+c^+}{4}}e^{C(1+\alpha^2)t},
  \end{align*}
  (note that \cite[Theorem 6.2]{Ouhabaz} apply also to non-symmetric forms, see \cite[Remark (ii), page 159]{Ouhabaz}). By duality the latter inequality imply also
  \begin{align*}
  	\|e^{t\mathcal B_{\alpha,\psi}}\|_{L^1_{\nu}\to L^2_\nu}\leq Ct^{-\frac{N+1+c^+}{4}}e^{C(1+\alpha^2)t},\quad \|e^{t\mathcal B^*_{\alpha,\psi}}\|_{L^1_{\nu}\to L^2_\nu}\leq Ct^{-\frac{N+1+c^+}{4}}e^{C(1+\alpha^2)t}.
  \end{align*}
Then combining  the latter estimates we obtain, for some possibly different constants $C>0$, 
  \begin{align*}
  	\|e^{t\mathcal B_{\alpha,\psi}}\|_{L^1_{\nu}\to L^\infty}\leq Ct^{-\frac{N+1+c^+}{2}}e^{C(1+\alpha^2)t},\quad \|e^{t\mathcal B^*_{\alpha,\psi}}\|_{L^1_{\nu}\to L^\infty}\leq Ct^{-\frac{N+1+c^+}{2}}e^{C(1+\alpha^2)t}.
  \end{align*}
The
 Dunford-Pettis criterion  yields then the  existence of a positive kernel $p_{\alpha,\psi,\nu}$ such that
\begin{align*}
	e^{t\mathcal B_{\alpha\psi}}f(z_1)=
\int_{\R^{N+1}_+}p_{\alpha,\psi,\nu}(t,z_1,z_2)f(z_2)\,d\nu,\quad f\in L^2_{\nu}
\end{align*}
and $0 \le p_{\alpha,\psi,\nu}(t,z_1,z_2) \le Ct^{-\frac{N+1+c^+}{2}}e^{C(1+\alpha^2)t}$ with $C>0$  independent of $\alpha, \psi$. Using the equivalence $e^{t\mathcal B_{\alpha, \psi}}=e^{-\alpha\psi}e^{t\tilde{\mathcal L}} e^{\alpha\psi}$ of equation \eqref{relation tilde L pert}, 
the heat kernel $p_{\tilde{\mathcal L}}$ of $\tilde{\mathcal L}$  satisfies
\begin{align*}
	0\leq 
e^{-\alpha \psi(z_1)}p_{\tilde{\mathcal L}}(t,z_1,z_2)e^{\alpha \psi(z_2)}=p_{\alpha,\psi,\nu}(t,z_1,z_2)
\leq 
Ct^{-\frac{N+1+c^+}{2}}e^{C(1+\alpha^2)t}.
\end{align*} 
Optimizing  on $\alpha \in \R$ we get 
\begin{align*}
0\leq 
p_{\tilde{\mathcal L}}(t,z_1,z_2)
\leq 
Ct^{-\frac{N+1+c^+}{2}}e^{Ct}\exp{\left(-\frac{(\psi(z_1)-\psi(z_2)^2}{4Ct}\right)}.
\end{align*}

Optimizing finally on $\psi \in \Psi$ and recalling Remark \ref{equiv distance} we get the required assertion. 
\qed

We can  now prove heat kernel estimates for the   degenerate   operator 
\begin{equation*} 
	\mathcal L=\Delta_x +2a\cdot \nabla_xD_y+ B_y
\end{equation*}
in $L^2_c$ with Neumann boundary conditions at $y=0$: $\lim_{y \to 0} y^c D_yu=0$.
 We recall that from Definition \ref{defi tilde L} we have  $\mathcal L=\phi \tilde {\mathcal L}\phi^{-1}$
and  then for $z_1=(x_1,y_1),z_1=(x_2,y_2)\in\R^{N+1}_+$
\begin{align*}
	e^{t\mathcal{L}}f(z_1)
	&=\phi(y_1) \big[e^{t\tilde{\mathcal L}}(\phi^{-1}f)\Big](z_1)
	=\int_{\R^{N+1}_+}\phi(y_1) p_{\tilde{\mathcal L}}(t,z_1,z_2)\phi(y_2)^{-1}f(z_2)\,d\nu(z_2)
	\\
	&:=\int_{\R^{N+1}_+}p_{\mathcal L}(t,z_1,z_2)f(z_2)\,y_2^cdz_2
\end{align*}
where 
\begin{align}\label{equiv kernel}
	p_{\mathcal L}(t,z_1,z_2)=\phi(y_1)p_{\tilde{\mathcal L}}(t,z_1,z_2)\phi(y_2).
\end{align}
Using the scaling property of $e^{t\mathcal L}$, see \eqref{scale.semi}, we obtain the following kernel estimates which constitute the main result of this section. Note that the distinction between the cases $c<0$ and $c\geq 0$ disappears. 

\begin{teo}\label{true.kernel} 
	Let $c+1>0$, $|a|<1$. The semigroup $(e^{t{\mathcal L}})_{t> 0}$ consists of integral operators i.e. there exists $p_{{\mathcal L}}(t,\cdot,\cdot)\in L^\infty(\R^{N+1}_+\times\R^{N+1}_+)$ such that  for $t>0$, $z_1=(x_1,y_1),\ z_2=(x_2,y_2)\in\R^{N+1}_+$
	\begin{align*}
		e^{t{\mathcal  L}}f(z_1)=
		\int_{\R^{N+1}_+}p_{{\mathcal L}}(t,z_1,z_2)f(z_2)\,y_2^c\,dz_2,\quad f\in L^2_{c}.
	\end{align*}
	Its heat kernel $p_{{\mathcal L}}$, written  with respect the measure $y^cdz$, satisfies for some  $C,k>0$, 
	\begin{align*}
		0\leq p_{{\mathcal  L}}(t,z_1,z_2)
		\leq C t^{-\frac{N+1}{2}} y_1^{-\frac{c}{2}} \left(1\wedge \frac {y_1}{\sqrt t}\right)^{\frac{c}{2}} y_2^{-\frac{c}{2}} \left(1\wedge \frac{y_2}{\sqrt t}\right)^{\frac{c}{2}}\,\exp\left(-\dfrac{|z_1-z_2|^2}{kt}\right).
	\end{align*}
The last estimate can be written  equivalently as
\begin{align*}
	0\leq p_{{\mathcal  L}}(t,z_1,z_2)
	\leq C t^{-\frac{N+1}{2}}  y_2^{-c} \left(1\wedge \frac{y_2}{\sqrt t}\right)^{c}\,\exp\left(-\dfrac{|z_1-z_2|^2}{kt}\right).
\end{align*}
\end{teo} 
{\sc Proof.} The existence of the heat kernel follows from Lemma \ref{realheat tilde} and  \eqref{equiv kernel}. To prove the required estimate we use  \eqref{Stima ptilde reale} for $t=1$ and    \eqref{Prop-phi} and  \eqref{equiv kernel} to get  
\begin{align*}
	p_{\mathcal{L}}(1,z_1,z_2)&\leq 
	C
	y_1^{-\frac{c}{2}}(1\wedge y_1)^{\frac{c}{2}}y_2^{-\frac{c}{2}}(1\wedge y_2)^{\frac{c}{2}}y_2^c\,
	\exp\left(\frac{-|z_1-z_2|^2}{k}\right)
\end{align*}
for any  $(z_1,z_2)\in \R^{N+1}_+\times \R^{N+1}_+.$
On the other hand using  the scaling property  \eqref{scale.semi} of Proposition \ref{Prop scaling} we get 
\begin{align*}
p_{\mathcal{L}}(s^2t,z_1,z_2)=s^{-(N+1+c)}p_{\mathcal{L}}(t,s^{-1}z_1,s^{-1}{z_2}),\qquad s,t>0, \quad {z_1}, {z_2}\in \R^{N+1}_+. 
\end{align*}
 In particular 
\begin{align*}
	p_{\mathcal{L}}(t,z_1,z_2)=t^{-\frac{N+1+c}2}p_{\mathcal{L}}\left(1,\frac{z_1}{\sqrt t},\frac{z_2}{\sqrt t}\right)
\end{align*}
which combined with the above inequality for $t=1$ gives the first required estimate.

The second equivalent estimate, up to a possibly different constants $C, k$, follows by applying \cite[Lemma 10.2]{MNS-Caffarelli}.
\qed

	\begin{os}
We remark that the above estimates seem to be optimal since in the case $N=0$ and $a=0$  coincide with the optimal estimate proved in \cite{MNS-Sharp}. Moreover, in the case $\mathcal{L}=\Delta_x+B_y$,  our upper bounds coincide with those of \cite{MNS-Caffarelli} which are  optimal.
\end{os}

\begin{os}\label{Oss c>0}
When $c\geq 0$, with the same proof of Lemma \ref{realheat tilde}, one could have obtained directly kernel estimates for $e^{t\mathcal L}$ (without relying on the weight $\phi$ and on the equivalent operator $\tilde{\mathcal L}$) by using  the Gagliardo-Nirenberg  inequalities of Proposition \ref{Sobolev3b} in place of Corollary \ref{Gagl-Nir-bpert}. This approach leads to the Gaussian estimate
\begin{align*}
	0\leq p_{{\mathcal  L}}(t,z_1,z_2)
	\leq C t^{-\frac{N+1+c}{2}}\exp\left(-\dfrac{|z_1-z_2|^2}{kt}\right)
\end{align*} 
which is weaker then the estimate proved in Theorem \ref{true.kernel}. 
\end{os}

With the same proofs of Lemma \ref{realheat tilde} and of Theorem \ref{true.kernel} we can  deduce heat kernel estimates for the degenerate operator   
\begin{equation*} 
	\mathcal L^\ast=\Delta_x +2a\cdot \nabla_xD_y+2c\frac {a\cdot \nabla_x u}{y}+ B_y
\end{equation*}
in $L^2_c$  with the oblique condition at $y=0$: $\lim_{y \to 0} y^c \left ( D_yu+2a\cdot \nabla_x u\right )=0$. 
\begin{teo}\label{real kernel dual} 
	Let $c+1>0$, $|a|<1$. The semigroup $(e^{t{\mathcal L^*}})_{t> 0}$ consists of integral operators i.e. there exists $p_{{\mathcal L^*}}(t,\cdot,\cdot)\in L^\infty(\R^{N+1}_+\times\R^{N+1}_+)$ such that  for $t>0$, $z_1=(x_1,y_1),\ z_2=(x_2,y_2)\in\R^{N+1}_+$
	\begin{align*}
		e^{t{\mathcal  L^\ast}}f(z_1)=
		\int_{\R^{N+1}_+}p_{{\mathcal L^\ast}}(t,z_1,z_2)f(z_2)\,y_2^c\,dz_2,\quad f\in L^2_{c}.
	\end{align*}
	Its heat kernel $p_{{\mathcal L}^\ast}$, written  with respect the measure $y^cdz$, satisfies for some  $C,k>0$, 
	\begin{align*}
		0\leq p_{{\mathcal  L}^\ast}(t,z_1,z_2)
		\leq C t^{-\frac{N+1}{2}} y_1^{-\frac{c}{2}} \left(1\wedge \frac {y_1}{\sqrt t}\right)^{\frac{c}{2}} y_2^{-\frac{c}{2}} \left(1\wedge \frac{y_2}{\sqrt t}\right)^{\frac{c}{2}}\,\exp\left(-\dfrac{|z_1-z_2|^2}{kt}\right).
	\end{align*}
	The last estimate can be written  equivalently as
	\begin{align*}
		0\leq p_{{\mathcal  L}^\ast}(t,z_1,z_2)
		\leq C t^{-\frac{N+1}{2}}  y_2^{-c} \left(1\wedge \frac{y_2}{\sqrt t}\right)^{c}\,\exp\left(-\dfrac{|z_1-z_2|^2}{kt}\right).
	\end{align*}
\end{teo}

Once kernel estimates have obtained for positive $t$, based on the argument of \cite[Theorem 6.16]{Ouhabaz} relying on the Phragmen-Lindel\"of theorem, we can extend them for complex times $t\in \Sigma_{\frac{\pi}{2}-\arctan\frac{|a|}{1-|a|}}$. We do not repeat the proof which is identical to that of  \cite[Section 5.3]{MNS-Singular-Half-Space}. For simplicity we state the result only for the operator $\mathcal{L}$.

 \begin{teo}\label{complexGauss} 
 	Let $c+1>0$, $|a|<1$. The semigroup $(e^{t{\mathcal L}})$, $t\in\Sigma_{\frac{\pi}{2}-\arctan\frac{|a|}{1-|a|}}$ consists of integral operators i.e. there exists $p_{{\mathcal L}}(t,\cdot,\cdot)\in L^\infty(\R^{N+1}_+\times\R^{N+1}_+)$ such that  for $t\in\Sigma_{\frac{\pi}{2}-\arctan\frac{|a|}{1-|a|}}$, $z_1=(x_1,y_1),\ z_2=(x_2,y_2)\in\R^{N+1}_+$
 	\begin{align*}
 		e^{t{\mathcal  L}}f(z_1)=
 		\int_{\R^{N+1}_+}p_{{\mathcal L}}(t,z_1,z_2)f(z_2)\,y_2^c\,dz_2,\quad f\in L^2_{c}.
 	\end{align*}
 Moreover for every $\eps>0$  there exist  constants $C_{\eps}, k_\eps>0$ such that if  $t\in\Sigma_{\frac{\pi}{2}-\arctan\frac{|a|}{1-|a|}-\eps}$ then
 	\begin{align*}
 		|p_{{\mathcal  L}}(t,z_1,z_2)|
 		\leq C_\epsilon |t|^{-\frac{N+1}{2}}  y_2^{-c} \left(1\wedge \frac{y_2}{\sqrt{ |t|}}\right)^{c}\,\exp\left(-\dfrac{|z_1-z_2|^2}{k_\epsilon|t|}\right).
 	\end{align*}
 \end{teo}

\section{General operators and oblique derivative}\label{Section oblique}
Kernel estimates for more general operators and boundary conditions follow by a linear change of variables, as we explain below. Let us first remove the assumption on the special form of  $\mathcal L=\Delta_{x} +2a\cdot\nabla_xD_y+ B_yu$ by considering the general operator in  $\R^{N+1}_+$

\begin{equation*}
	\mathcal L =\mbox{Tr }\left(QD^2u\right)+\frac{c}{ y}D_y=\sum_{i,j=1}^{N}q_{ij}D_{x_ix_j}+2\sum_{i=1}^Nq_iD_{x_iy}+\gamma D_{yy}+\frac{c}{y}D_y
\end{equation*}
with Neumann boundary condition $\lim_{y \to 0} y^{\frac c\gamma} D_yu=0$. 
If $Q_1$ is the $N\times N$ matrix $(q_{ij})$, $\gamma=q_{N+1,N+1}$ and $q=(q_1, \dots, q_N)$, we assume that the quadratic form $Q(\xi,\eta)=Q_1(\xi,\xi)+\gamma \eta^2+2 q\cdot\xi \,\eta$ is positive definite.
Through a linear change of variables in the $x$ variables the term $\sum_{i,j=1}^{N}q_{ij}D_{x_ix_j}$ is transformed into $\gamma\Delta_x$ and all the results of the previous sections remain true. In particular if $\frac c \gamma+1>0$ then  $\mathcal L$  generates a contractive analytic semigroups $e^{z \mathcal L}$
in $L^2_{\frac c \gamma}$ and  the kernel estimates of Theorems \ref{true.kernel} and \ref{complexGauss}   remain true.

A further change of variables allows to deal with  the operator 
\begin{equation*}
	\mathcal L =\mbox{Tr }\left(QD^2u\right)+\frac{ v\cdot \nabla }{y}=\sum_{i,j=1}^{N+1}q_{ij}D_{ij}+\frac{b\cdot \nabla_x}{y}+\frac{cD_y}{y},\qquad c\neq 0, 
\end{equation*}
where $v=(b,c)\in\R^{N+1}$ and is $Q$ positive definite. We impose, in this case,  an oblique derivative boundary condition \begin{align*}
	\lim_{y \to 0} y^{\frac c \gamma}\, v \cdot \nabla u=0.
\end{align*} 
We transform $\mathcal L$ into a similar operator with $b=0$ and Neumann boundary conditions by defining the following isometry of $L^2_{\frac c\gamma}$
\begin{align}\label{Tran def}
	T\, u(x,y)&:=u\left(x-\frac b c y,y\right),\quad (x,y)\in\R_+^{N+1}.
\end{align}

\begin{lem}\label{Isometry action der} Let $1< p< \infty$,  $ v=(b,c)\in\R^{N+1}$,  $c\neq 0$. Then  for $u\in W^{2,1}_{loc}\left(\R_+^{N+1}\right)$ 
\begin{itemize}
	\item[(i)] $\ds T^{-1}\,\left(v\cdot \nabla \right) T\,u=cD_yu$;
	
	\item[(ii)] $\ds 	T^{-1}\,\left(\mbox{Tr }\left(QD^2\right)+\frac{ v\cdot \nabla  }{y}\right) T\,u=\mbox{Tr }\left(\tilde QD^2u\right)+\frac{c}{y}D_y u$.
\end{itemize}
Here $\tilde Q$  is a uniformly elliptic symmetric matrix defined by 
		\begin{align*}
			\tilde Q=
			\left(
			\begin{array}{c|c}
				Q_N-\frac{2}c b\otimes q{+\frac \gamma{c^2} b\otimes b}  & q^t- \frac{\gamma}c b^t \\[1ex] \hline
				q- \frac{\gamma}c b & \gamma
			\end{array}\right)
		\end{align*} 
		and  $\gamma= q_{N+1,N+1}$.
\end{lem}{\sc{Proof.}} The proof follows by a straightforward computation.
\qed
\medskip

The latter lemma shows that the operator $\mathcal L$, endowed with the oblique  boundary condition $\lim_{y \to 0} y^{\frac c \gamma}\, v \cdot \nabla u=0$ is unitarily equivalent to  the operator
\begin{align}\label{equiv oblique}
		 \mathcal {\tilde{L}}:=T^{-1}\mathcal L T\,=\mbox{Tr }\left(\tilde QD^2\right)+\frac{c}{y}D_y 
	\end{align}
with the Neumann  boundary condition $\lim_{y \to 0} y^{\frac c \gamma}\, D_yu=0$.  Then, by similarity and the previous results,  if $\frac c \gamma+1>0$,    $\mathcal L$  generates a contractive analytic semigroups $e^{z \mathcal L}$. In order to deduce kernel estimates we observe that, if $p_{\tilde{\mathcal L}}$ is the heat kernel of $\tilde{\mathcal L}$ with respect to the measure $y^{\frac c\gamma}dz$, then one has for $z_1=(x_1,y_1),z_2=(x_2,y_2)\in\R^{N+1}_+$  
\begin{align*}
	e^{t\mathcal{L}}f(z_1)
	&=T \big[e^{t\tilde{\mathcal L}}(T^{-1}f)\Big](z_1)=\big[e^{t\tilde{\mathcal L}}(T^{-1}f)\Big](x_1-\frac{b}{c}y_1, y_1)\\[1ex]
	&	=\int_{\R^{N+1}_+} p_{\tilde{\mathcal L}}\left(t,\left(x_1-\frac{b}{c}y_1,y_1\right),z_2\right)f\left(x_2+\frac{b}{c}y_2,y_2\right)\,y_2^{\frac c \gamma }dx_2dy_2
	\\[1ex]
	&=\int_{\R^{N+1}_+}p_{\tilde{\mathcal L}}\left(t,\left(x_1-\frac{b}{c}y_1,y_1\right),\left(x_2-\frac{b}{c}y_2, y_2\right)\right)f(z_2)\,y_2^{\frac c \gamma }dz_2.
\end{align*}
This shows that  the heat kernel of $\mathcal L$ is given by 
\begin{align}\label{equiv kernel}
	p_{\mathcal L}(t,z_1,z_2):=p_{\tilde{\mathcal L}}\left(t,\left(x_1-\frac{b}{c}y_1,y_1\right),\left(x_2-\frac{b}{c}y_2, y_2\right)\right).
\end{align}
By Theorem \ref{true.kernel}, we then obtain kernel estimates for $\mathcal L$. We state them for real $t$ only for simplicity. We need a preliminary lemma.
\begin{lem} \label{gaussian-exp}
	For any $k\in\R^N\setminus\{0\}$ one has, for some positive constant $C>0$,
	\begin{align*}
		|x-ky|^2+y^2\geq C\left(|x|^2+y^2\right),\qquad \forall x\in\R^n,\, y\in\R.
	\end{align*}
\end{lem}
{\sc{Proof.}} Using the  Cauchy-Schwarz inequality we get for $\epsilon>0$
\begin{align*}
	|x-ky|^2+y^2&=|x|^2+(1+|k|^2)y^2-2x\cdot k y\geq|x|^2+(1+|k|^2)y^2-2|x||k| |y|\\[1ex]
	&\geq|x|^2+(1+|k|^2)y^2-\frac 1 \epsilon |x|^2-\epsilon {|k|^2} y^2\\[1ex]
	&=(1-\frac{1}\epsilon)|x|^2+(1+(1-\epsilon)|k|^2)y^2.
\end{align*}
The required claim then follows by choosing $1<\epsilon<1+\frac{1}{|k|^2}$. 
\qed

\begin{teo} \label{complete-oblique}
	Let $v=(b,c)\in\R^{N+1}$ with $b=0$ if $c =0$, $
	\left(
	\begin{array}{c|c}
		Q  & { q}^t \\[1ex] \hline
		q& \gamma
	\end{array}\right)$
	an  elliptic matrix. 
If $\frac{c}{\gamma}+1>0$ then the operator
	\begin{align*}
		\mathcal L =\mbox{Tr }\left(QD^2u\right)+\frac{ v\cdot \nabla }{y}
	\end{align*} endowed with the boundary conditions
\begin{align*}
	&\lim_{y \to 0} y^{\frac c \gamma}\, v \cdot \nabla u=0,\qquad \text{(if $c\neq 0$)},\qquad \qquad 
	\lim_{y \to 0} y^{\frac c \gamma}\, D_y u=0,\qquad \hspace{2ex}\text{(if $c=0$)},
\end{align*}
 generates a contractive analytic semigroup  in $L^2_{\frac c\gamma}$. The semigroup $(e^{t{\mathcal L}})_{t> 0}$ consists of integral operators	and its heat kernel $p_{{\mathcal L}}$, written  with respect the measure $y^\frac{c}{\gamma}dz$, satisfies for some  $C,k>0$, 
 \begin{align*}
 	0\leq p_{{\mathcal  L}}(t,z_1,z_2)
 	\leq C t^{-\frac{N+1}{2}} y_1^{-\frac{c}{2\gamma}} \left(1\wedge \frac {y_1}{\sqrt t}\right)^{\frac{c}{2\gamma}} y_2^{-\frac{c}{2\gamma}} \left(1\wedge \frac{y_2}{\sqrt t}\right)^{\frac{c}{2\gamma}}\,\exp\left(-\dfrac{|z_1-z_2|^2}{kt}\right).
 \end{align*}
\end{teo}
{\sc{Proof.}} 	The generation property follows by similarity using the equality \eqref{equiv oblique}. To prove the kernel estimates we use \eqref{equiv kernel} and Theorem \ref{true.kernel} getting 
\begin{align*}
	0\leq p_{{\mathcal  L}}(t,z_1,z_2)
	\leq C t^{-\frac{N+1}{2}} y_1^{-\frac{c}{2\gamma}} &\left(1\wedge \frac {y_1}{\sqrt t}\right)^{\frac{c}{2\gamma}} y_2^{-\frac{c}{2\gamma}} \left(1\wedge \frac{y_2}{\sqrt t}\right)^{\frac{c}{2\gamma}}
	\\[1ex]
	&\times\exp\left(-\dfrac{|x_1-x_2-\frac{b}{c}(y_1-y_2)|^2+|y_1-y_2|^2}{kt}\right).
\end{align*}
The required claim then follows using  Lemma \ref{gaussian-exp}.
\qed
\medskip 
Having at disposal good kernel estimates for $e^{t\mathcal L}$ implies various consequences about the regularity of the semigroup and its spectral properties; as an application   one can study, for example, generation  results in $L^p$ space. 
 We refer the reader to \cite{MNS-Singular-Half-Space, MNS-Degenerate-Half-Space} where it has been proved that  the above semigroup extrapolates to the weighted spaces $L^p\left(\R^{N+1}_+, y^mdxdy\right)$, $m\in\R$, when   $0<\frac{m+1}p<\frac{c}\gamma+1$; furthermore a characterization of the domain of the generators is given in terms of  suitable anisotropic weighted Sobolev spaces.

\bibliography{../../TexBibliografiaUnica/References}
\end{document}